\documentclass[10pt]{report}
\usepackage[utf8]{inputenc}
\usepackage{latexsym,makeidx,amssymb,amsmath,theorem,url}
\usepackage[makeroom]{cancel}
\usepackage{imakeidx}
\usepackage{yfonts}

\newtheorem{lemma}{Lemma}[section]
\newtheorem{thm}[lemma]{Theorem}
\newtheorem{cor}[lemma]{Corollary}
\newtheorem{prop}[lemma]{Proposition}
\newtheorem{defi}[lemma]{Definition}

{\theorembodyfont{\rmfamily}}
\newcommand{\N}{\mathbb{N}}

\newcommand{\Z}{\mathbb{Z}}
\newcommand{\Q}{\mathbb{Q}}
\newcommand{\R}{\mathbb{R}}

\newcommand{\C}{\mathbb{C}}
\newcommand{\e}{\mathrm{e}}

\newcommand{\ord}{\mathrm{ord}}

\newcommand{\ep}{\varepsilon}

\newcommand{\al}{\alpha}
\newcommand{\be}{\beta}

\newcommand{\duk}{\par\noindent{\bf Proof. }}
\newcommand{\kduk}{\hfill{$\Box$}\bigskip}

\newcommand{\sus}{\subset}

\newcommand{\cc }{\colon}

\renewcommand\contentsname{Contents\\
{\small{\rm (${}^c\ds$ means a~relatively complete section.)}}}

\def\ds{\dots}

\DeclareFontFamily{U}{wncy}{}
    \DeclareFontShape{U}{wncy}{m}{n}{<->wncyr10}{}
    \DeclareSymbolFont{mcy}{U}{wncy}{m}{n}
    \DeclareMathSymbol{\Sh}{\mathord}{mcy}{"58} 

\begin{document}
\pagenumbering{roman}
\thispagestyle{empty}
\begin{center}
    {\Huge {\bf Catalan's conjecture is\\

    \medskip
    Mih\u ailescu's theorem}}
    
    \bigskip\bigskip
    {\huge Martin Klazar}

    \bigskip
    {\Large (KAM MFF UK Praha)}
    
\end{center}

\newpage
\thispagestyle{empty}
\noindent
dedicated to my parents Blanka and Ji\v r\'\i

\newpage
\thispagestyle{empty}
\noindent
\centerline{$-$ 520 $-$}

\medskip
{\bf 48. {\em Th\'eor\`eme}. Deux nombres entiers cons\'ecutifs, autres\\
que 8 et 9, ne peuvent \^etre des puissances exactes. (Catalan.)}

\medskip\noindent
See \cite{cata1}.

\bigskip\noindent
192

\centerline{{\bf 13.}}

\centerline{{\bf Note}}
\centerline{{\footnotesize extraite d'une lettre adress\'ee \`a l'\'editeur par Mr. {\em
E. Catalan}, R\'ep\'etiteur \`a l'\'ecole}}
\vspace{-1mm}
\centerline{{\footnotesize polytechnique de Paris.}}
\vspace{-1mm}
\centerline{\rule{1.8cm}{0.3pt}}

\smallskip\noindent
{\footnotesize {\bf ,,\,$\mathbb{J}$e vous prie, Monsieur, de vouloir bien
\'enconcer, dans votre recueil, le\\
,,\,th\'eor\`eme suivant, que je crois vrai, bien que je n'aie pas encore r\'eussi \`a\\
,,\,le d\'emontrer compl\`etement: d'autres seront peut-\^etre plus
heureux:

,,\,Deux nombres entiers cons\'ecutifs, autres
que 8 et 9, ne peuvent \^etre

,,\,des puissances exactes; autrement dit: l'\'equation $\mathbf{x^m-y^n=1}$, dans

,,\,laquelle les inconnues sont enti\`eres et positives, n'adm\`et qu'une seule

\vspace{-1mm}
,,\,solution.\,''}}

\vspace{-1mm}
\centerline{\rule{3cm}{0.5pt}}

\medskip\noindent
See \cite{cata2}.

\bigskip\bigskip\noindent
---\,written according to \cite[pp. 1--2]{bilu_al}. In 1842, the Belgian-French mathematician Eug\`ene Ch. 
Catalan (1814--1894) (see \cite{cataWiki})
put forth without a~proof a~theorem asserting that $8$ and $9$ 
are the only consecutive pure powers. Two years later he corrected himself  and changed it to a~conjecture. Catalan's conjecture was proven in 2004 by the Romanian 
mathematician Preda Mih\u ailescu (1955). In this text, we present 
a~complete proof of Catalan's conjecture. In particular, we give in entirety Mih\u{a}ilescu's proof.
\newpage

\section*{Introduction}
\addcontentsline{toc}{chapter}{Introduction}
One can state Catalan's conjecture \cite{cata2} (1844) as follows. For integers $m,n\ge2$ 
the only solution of the Diophantine equation 
$$
x^m-y^n=1
$$
in nonzero integers $x,y$ is $m=y=2$, $x=\pm3$ and $n=3$. After many partial results by many authors, Catalan's conjecture was proven in 2004 by 
P.~Mih\u ailescu \cite{miha}.

Three excellent books on
the subject are \cite{ribe} by P.~Ribenboim, \cite{scho} by R.~Schoof and 
\cite{bilu_al} by Yu. Bilu, Y.~Bugeaud and M.~Minotte. The first book is from 
the pre-Mih\u ailescu era; the other two present Mih\u ailescu's 
proof. They are far from self-contained, though. In this text our 
aim is to give a~complete and self-contained presentation of 
Mih\u ailescu's proof, and in fact of the proof of Catalan's conjecture (which has to include many things 
Mih\u ailescu took for granted). The initial part in Chapters~\ref{chap_euler}--\ref{sec_?} 
corresponds to my 
lectures in the course {\em Algebraic Number Theory} in 2024/25 and 2025/26.

It is easy to see that the exponents $m=p$ and $n=q$ can be assumed to be 
distinct primes. The resolution of Catalan's conjecture naturally divides 
in the elementary part with $p=2$ or $q=2$, and the non-elementary 
part with $p,q\ge3$. ``Elementary'' does not mean easy, and 
``non-elementary'' means hard. The elementary part is covered by the 
first three chapters, and the rest of our text is devoted to the 
non-elementary case. If it is not said else, a~solution is always an 
integral solution.

In Chapter~\ref{chap_euler}, we resolve the equation $x^2-y^3=1$; the only solutions are $\langle\pm3,2\rangle$, $\langle\pm1,0\rangle$, and 
$\langle0,-1\rangle$. We present three resolutions of the equation. The historically 
first due to L.~Euler (1737) is in Section~\ref{sec_descent}. It makes 
use of the fact that, in the modern view, $x^2=y^3+1$ is an elliptic 
curve. Euler actually proved that the five mentioned solutions are the 
only rational solutions. The second resolution of $x^2-y^3=1$ in 
Section~\ref{sec_klaz} is due to this author in 1989; the factorization $x^2=
(y+1)(y^2-y+1)$ leads to investigation of properties of solutions of the Pell 
equation $x^2-3y^2=1$. Another way how to solve $x^2-y^3=1$ is to start with 
the factorization $(x+1)(x-1)=y^3$; one is then led to equations $x^3-
2y^3=\pm1$. In 1957, in a~little known article, A. Wakulicz provided an elementary 
resolution of the more general equation $x^3+y^3=2z^3$; we 
present his result in Section \ref{sec_wakulicz}. All three 
resolutions of $x^2-y^3=1$ are completely elementary in the sense 
that they take place in the field $\Q$.  

% Melody, the music, of proofs.

\bigskip\bigskip\noindent
Praha and Louny, December 2025 to ??\hfill{Martin Klazar}

\newpage

\section*{Notation}
\addcontentsline{toc}{chapter}{Notation}

$\N=\{1,2,\ds\}$ is the set of natural numbers and $\N_0=\{0,1,2,\ds\}$. 

\newpage
\tableofcontents
\newpage

\pagenumbering{arabic}
\chapter[Euler's theorem: $x^2-y^3=1$]{Euler's theorem}\label{chap_euler}

We begin our long journey on the mountain range of the proof of 
Catalan's conjecture in the first elementary case: the only integral solutions of the 
equation
$$
x^2-y^3=1
$$
are pairs $\langle\pm3,2\rangle$, $\langle\pm1,0\rangle$ and 
$\langle0,-1\rangle$. We give three proofs. Section~\ref{sec_pell} collects 
auxiliary results on Pell equations $x^2-dy^2=1$ needed in this chapter and in Chapter~\ref{chap_chaoKo}. In 
Section~\ref{sec_UFDaPS} we survey unique factorization domains and 
principles of powers. In Section~\ref{sec_klaz} we present the resolution of $x^2-y^3=1$ in \cite{klaz}, found by this author in 1989. Section~\ref{sec_descent}
is devoted to a~modern proof of the result due to Leonhard Euler 
(1707--1783) (\cite{eulerWiki}) in \cite{eule} in 1738 that the five pairs are 
the only {\em rational} solutions  of the equation. In the last Section~\ref{sec_wakulicz} we present 
the elementary proof of A.~Wakulicz \cite{waku} that the equation 
$x^3+y^3=2z^3$ has no solution with $x\ne\pm y$. Thus the only 
solutions of $x^3-2y^3=\pm1$ are $\langle\pm1,0\rangle$ and 
$\langle\pm1,\pm1\rangle$ (equal signs). Using the factorization
$$
(x-1)(x+1)=y^3
$$
one deduces that $x^2-y^3=1$ has just 
the five mentioned solutions. The existence of the nontrivial (nonzero) 
solution
$\langle\pm3,2\rangle$ gives
this elementary case of Catalan's conjecture a distinct flavor compared to the elementary cases in the next two chapters. In 
Sections~\ref{sec_klaz}, \ref{sec_descent} and 
\ref{sec_wakulicz} we follow \cite{klaz}, \cite{conr} and 
\cite{waku}, respectively.

\section[${}^c$Pell equations]{Pell equations}\label{sec_pell}

{\em Pell equation} is any Diophantine equation of the form
$$
x^2-dy^2=1\,,
$$
with unknowns $x,y$ and parameter $d\in\N$ that is not a~square: $d=2$, 
$3$, $5$, $6$, $7$, $8$, $10$ and so on. We say that a~pair $a,b\in\N$ is 
a~{\em minimal solution} of the equation if $a^2-db^2=1$ and there is no solution $x,y\in\N$ with $x<a$. Minimal solutions, if they exist, are 
unique. 

\begin{prop}\label{prop_PellEq}
Let $a,b\in\N$ be the minimal solution of the Pell equation
$$
x^2-dy^2=1\,.
$$
Then natural solutions $x,y\in\N$ of the equation form an infinite set
$$
\{\langle x,\,y\rangle\in\N^2\cc\;\text{$x+y\sqrt{d}=\big(a+b\sqrt{d}\big)^n$ for some $n\in\N$}\}\,.
$$
\end{prop}
\duk
Let $x_i,y_i\in\Z$, $i=1,2$, be two solutions of the equation and $x_3,y_3\in\Z$ be defined by
$$
x_3+y_3\sqrt{d}=
\big(x_1+y_1\sqrt{d}\big)\big(x_2+y_2\sqrt{d}\big)\,.
$$
Then also
$$
x_3-y_3\sqrt{d}=
\big(x_1-y_1\sqrt{d}\big)\big(x_2-y_2\sqrt{d}\big)\,.
$$
Multiplying the displayed equalities we get
$$
x_3^2-dy_3^2=\big(x_1^2-dy_1^2\big)
\big(x_2^2-dy_2^2\big)=1\cdot1=1\,.
$$
Thus $x_3,y_3$ is a~solution of the equation as well. Even more easily, if $x,y$ is a~solution of the equation, then the reciprocal
$$
\frac{1}{x+y\sqrt{d}}=x-y\sqrt{d}
$$
produces solution $x,-y$.
Note that if $x_i,y_i\in\N$, $i=1,2$, are solutions of the equation, then
$$
x_1<x_2\iff x_1+y_1\sqrt{d}<x_2+y_2\sqrt{d}\,.
$$
Also, if $x,y\in\Z$ is a~solution of the equation, then $x,y\in\N$ iff $x+y\sqrt{d}>1$.

Now let $a,b$ be as stated and  $x,y\in\N$ be such that $x^2-dy^2=1$. We take the unique $m\in\N_0$ such that
$$
\al:=\big(a+b\sqrt{d}\big)^m<x+y\sqrt{d}\le\big(a+b\sqrt{d}\big)^{m+1}\,.
$$
If the last inequality were strict, then
$$
1<
u+v\sqrt{d}:=\al^{-1}\cdot(x+y\sqrt{d})<
a+b\sqrt{d}
$$
would give a~solution $u,v\in\N$ of the equation, in contradiction with the minimality of $a,b$.
Thus $x+y\sqrt{d}=\big(a+b\sqrt{d}\big)^{m+1}$.
\kduk

\noindent
In 1770, Lagrange proved that every Pell equation has a~minimal 
solution, and hence infinitely many solutions. Probably, 
we do not need this result to
resolve Catalan's conjecture. But we need some properties of the solutions of the Pell equation $x^2-3y^2=1$.

\begin{cor}\label{prop_x2_3y2}
The solutions of Pell equation
$$
x^2-3y^2=1
$$ 
are exactly the pairs $\langle\pm x_n,\pm y_n\rangle$, $n\in\N_0$, where
$$
x_n+y_n\sqrt{3}=\big(2+\sqrt{3}\big)^n\,.
$$
Also, $x_0=1$, $y_0=0$ and $x_{n+1}=2x_n+3y_n$, $y_{n+1}=x_n+2y_n$.
\end{cor}
\duk
This follows from Proposition~\ref{prop_PellEq} because $x^2-
3y^2=1$ has the minimal solution $\langle2,1\rangle$.
\kduk
\begin{prop}\label{prop_x2_3y2Q}
If $\langle x_n,y_n\rangle$, $n\in\N_0$, are as in the previous 
corollary, then for every $n\in\N_0$ we have
$$
x_{2n}=2x_n^2-1,\ y_{2n}=2x_ny_n,\ 
x_{2n+1}=(y_n+y_{n+1})^2+1,\ y_{2n+1}=2x_ny_{n+1}-1\,.
$$  
The numbers $x_n$ and $y_n$ have different parity, and $x_n$ is odd iff $n$ is even.
\end{prop}
\duk
Since $x_n^2-3y_n^2=1$, 
\begin{eqnarray*}
x_{2n}+y_{2n}\sqrt{3}&=&
\big(2+\sqrt{3}\big)^{2n}=\big(x_n+y_n\sqrt{3}\big)^2\\
&=&x_n^2+3y_n^2+2x_ny_n\sqrt{3}=
2x_n^2-1+2x_ny_n\sqrt{3}\,.
\end{eqnarray*}
Similarly,
\begin{eqnarray*}
x_{2n+1}+y_{2n+1}\sqrt{3}&=&(2+\sqrt{3})(x_n^2+3y_n^2+2x_ny_n\sqrt{3})\\
&=&2x_n^2+6x_ny_n+6y_n^2+(x_n^2+4x_ny_n+3y_n^2)\sqrt{3}\,.
\end{eqnarray*}
Now 
$2x_n^2+6x_ny_n+6y_n^2=x_n^2+6x_ny_n+9y_n^2+1=(y_n+x_n+2y_n)^2+1=
(y_n+y_{n+1})^2+1$ and $x_n^2+4y_ny_n+3y_n^2=2x_n^2+4x_ny_n-1=2x_n(x_n+2y_n)-1=2x_ny_{n+1}-1$. The last claim is immediate from the just proven formulas.
\kduk

\section[${}^c$UFD and PP]{UFD and PP}\label{sec_UFDaPS}

These acronyms refer to {\em unique factorization domain(s)} and {\em 
principle(s) of powers}, respectively. PP propel resolutions of Diophantine 
equations, and therefore deserve more 
than the usual glossing over. Here we treat them in detail.

\begin{prop}[PP0]\label{prop_aDelBc}
Let $a,b,c\in\Z$ and $k\in\N$. If $a$ divides $bc^k$ and $a,c$ are coprime, then $a$ divides $b$.     
\end{prop}

\begin{prop}[PP1]\label{prop_PP1}
Let $k,l\in\N$ with $k,l\ge2$. If $a_i,b\in\N_0$ for $i\in[k]$ are  
$k+1$ numbers such that the $a_i$ are pairwise coprime and if
$$
a_1a_2\ds a_k=b^l\,,
$$
then there exist $k$ pairwise coprime numbers $b_i\in\N_0$ such that for every $i\in[k]$, 
$$
a_i=b_i^l\,.
$$
If $l$ is odd then this result holds also when $\N_0$ is replaced with $\Z$.
\end{prop}

\begin{prop}[PP2]\label{prop_PP2}
Let $p$ be a~prime and $k\in\N$ with $k\ge2$. If $a,b,c$ in $\N_0$ are numbers such that 
$$\mathrm{gcd}(a,b)=p\wedge
ab=c^k\,, 
$$
then there exist coprime numbers $d,e\in\N_0$ such that
$$
\{a,\,b\}=\{pd^k,\,p^{k-1}e^k\}\,.
$$
If $k$ is odd then this result holds also when $\N_0$ is replaced with $\Z$.
\end{prop}

In order to generalize PP0, PP1 and PP2 to domains and to prove them, we review the notion of 
a~unique factorization domain. We also introduce irreducible 
factorizations. Recall that an {\em (integral) domain}
$$
R=\langle R,\,0_R,\,1_R,\,+,\,\cdot\rangle
$$
is a~commutative ring with $1_R$ such that for every $a,b\in R^*$ ($=R\setminus
\{0_R\}$) we have $ab\ne0_R$. For $a,b\in R$ we say that {\em $a$ divides $b$ (in $R$)}, written $a\,|\,b$, 
if $b=ac$ ($=a\cdot c$) for some $c\in R$. We say that 
$a\in R$ is a~{\em unit} if $a\,|\,1_R$, that is, $a$ is 
multiplicatively invertible. The set of units 
in $R$ is denoted by $R^{\times}$. It is easy to see that
$$
\langle R^{\times},\,1_R,\,\cdot\rangle
$$
is an Abelian group, the {\em group of units} of the domain $R$. For 
$a,b\in R$ we write $a\sim b$ if $a=bc$ for some $c\in R^{\times}$; we
say that the elements $a$ and $b$ are {\em associated}.
For example, in the domain of integers 
$$
\Z=\langle\Z,\,0,\,1,\,+,\,\cdot\rangle
$$
we have $m\sim n$ iff $m=\pm n$.
It is easy to see that $\sim$ is an equivalence relation, and that it is congruent with respect to multiplication. For $a\in R$
we denote by 
$$
[a]_{\sim}\ \ (=\{b\in R\cc\;b\sim a\})
$$
the {\em block} of the element $a$ in the equivalence $\sim$.
We get the (commutative) monoid of blocks
$$
\langle R/\!\sim,\,[1_R]_{\sim}=R^{\times},\,\cdot\rangle\,.
$$

\begin{prop}\label{prop_mutuDiv}
In any domain $R$, two elements are associated if and only if each divides the other.    
\end{prop}
\duk
Let $a,b\in R$. Suppose that $a\sim b$. Thus $a=bc$ and $ac^{-1}=b$ for some $c,c^{-1}\in R^{\times}$. Hence $b\,|\,a$ and $a\,|\,b$.

Suppose that $a\,|\,b$ and $b\,|\,a$. Thus $b=ac$ and $a=bd$ for some $c,d\in R$. We get the equality
$$
b\cdot(1_R-d\cdot c)=0_R\,.
$$
Thus, since we are in a~domain, $b=a=0_R$ and $a\sim b$, or $dc=1_R$ and again $a\sim b$.
\kduk

Two elements $a,b\in R$ are
{\em coprime}, written $(a,b)=1_R$, if they can be simultaneously divided 
only by units. 

\begin{defi}[gcd]\label{def_gcd}
Let $R$ be a~domain and $a,b\in R$. We say that $c\in R$ is 
the greatest common divisor of $a$ and $b$, and write $c=\mathrm{gcd}(a,b)$, 
if $c$ divides $a$ and $b$, and every simultaneous divisor of $a$
and $b$ divides $c$.   
\end{defi}
If $c=\mathrm{gcd}(a,b)$ and $c'\sim c$, then $c'=\mathrm{gcd}(a,b)$. If $a\sim a'$, $b\sim b'$,
$c=\mathrm{gcd}(a,b)$ and 
$c'=\mathrm{gcd}(a',b')$, then $c\sim c'$. Also, $a$ and $b$ are coprime iff $\mathrm{gcd}(a,b)=1_R$.

Let $R$ be a~domain.
An element $a\in R$ is 
{\em irreducible} if $a\in R^*\setminus R^{\times}$ and if in every
multiplicative decomposition 
$a=bc$ with $b,c\in R$,
$b$ or $c$ is a~unit. If $a,b\in R$, $a\sim b$ and $a$ is irreducible, then so is $b$. We denote the set of 
irreducibles in $R$ by
$R^{\mathrm{ir}}$.

\begin{defi}[UFD~1]\label{def_UFD1}
$R$ is a~unique factorization domain, or 
{\em UFD}, if every element in  $R^*\setminus R^{\times}$ is a~product of irreducibles, and this product is unique up to the order of factors and the relation $\sim$.  
\end{defi} 
In more details, $R$ is UFD if for every element $a\in R^*\setminus R^{\times}$ there exist $m\in\N$ irreducibles $a_i$ such that 
$$
a=a_1\cdot a_2\cdot\ldots\cdot a_m\,,
$$
and if any equality
$$
c_1\cdot c_2\cdot\ldots\cdot c_l= b_1\cdot b_2\cdot\ldots\cdot b_m\,,
$$
where $l,m\in\N$ and $b_i$ and $c_i$ are irreducibles, implies that $l=m$ and that there exists a~permutation $\pi$ of the numbers $1,2,\dots,l$ such that for
every $i\in[l]$ we have
$c_i\sim b_{\pi(i)}$. A~prototypical example of UFD is $\Z$, 
which we prove in Section~\ref{sec_eucUFD}.

Recall that if $A$ and $B$ are sets and $X\sus A\times B$, then the relation $X$ is 
a~{\em partial function} (from $A$ to $B$) if for every $a\in A$ there exists at most one 
$b\in B$ such that $\langle a,b\rangle\in X$. We formalize (actually, set-theorize) irreducible factorizations.

\begin{defi}[irreducible factorizations]\label{def_irrFact}
Let $R$ be a~domain, 
$\mathbb{I}=R^{\mathrm{ir}}/\!\sim$ and let $a\in R^*$. Irreducible factorizations (in $R$) are 
finite partial functions $X$ from $\mathbb{I}$ to $\N$.
If $X\ne\emptyset$, we say that $X$ is an 
irreducible factorization of the element $a$ if 
$$
[a]_{\sim}=\prod_{\langle\al,\,m\rangle\in X}
\al^m\,.
$$
We say that $X=\emptyset$ is an irreducible factorization of $a$ if $a\in R^{\times}$.
\end{defi}

\noindent
For any irreducible factorization $X$ we set
$$
(X)_1=\{\al\in\mathbb{I}\cc\;\exists m\in\N:\,\langle\al,\,m\rangle\in X\}\,,
$$
and for $\al\in(X)_1$ we denote by $X(\al)$ the unique $m\in\N$ such that
$\langle\al,m\rangle\in X$.

Let $a,b\in R$ and $X,Y$ be the respective irreducible factorizations.
We define that {\em $X$ divides $Y$}, written $X\,|\,Y$, if 
$(X)_1\sus(Y)_1$ and for every $\al\in(X)_1$ we have $X(\al)\le 
Y(\al)$. It is easy to see that $a\,|\,b$ iff $X\,|\,Y$.
We restate the definition of UFD.

\begin{defi}[UFD~2]\label{def_UFD2}
A~domain $R$ is a~unique factorization domain, or 
{\em UFD}, if every element in $R^*$ has a~unique irreducible factorization.
\end{defi}

\begin{prop}\label{prop_UFDgcd}
In every {\em UFD} $R$ every two elements $a,b\in R^*$ have the greatest common divisor $c$.  If $a$ and $b$ are coprime then $c=1_R$. Else, denoting by $X$ and $Y$ the irreducible factorization of $a$ and $b$, respectively, we have 
$$
[c]_{\sim}=
\prod_{\al\in X_1\cap Y_1}\al^{\min(X(\al),\,Y(\al))}\,.
$$  
\end{prop}
\duk
This is immediate from the interpretation of divisibility of elements in $R^*$ in terms of their irreducible factorizations. 
\kduk

\noindent
As for the pairs $a,b$ with $ab=0_R$, we have $\mathrm{gcd}(0_R,a)=a$ for every $a\in R^*$, and 
$\mathrm{gcd}(0_R,0_R)$ does not exist.

We generalize PP0, PP1 and PP2 to UFD. Let $R$ be UFD, $a,b\in R^*$ and let $X$ and $Y$ be the respective irreducible factorizations of $a$ and 
$b$. Then $a$ and $b$ are coprime iff $(X)_1\cap(Y)_1=\emptyset$. 
The product $ab$ has the irreducible factorization 
\begin{eqnarray*}
XY&:=&\{\langle\al,
\,X(\al)\rangle\cc\;\al\in (X)_1\setminus(Y)_1\}\cup\{\langle\al,
\,Y(\al)\rangle\cc\;\al\in (Y)_1\setminus(X)_1\}\,\cup\\
&&\cup\,\{\langle\al,
\,X(\al)+Y(\al)\rangle\cc\;\al\in (X)_1\cap(Y)_1\}\,.
\end{eqnarray*}
If $a$ and $b$ are coprime then we have the disjoint union $XY=X\cup Y$. 
For $l\in\N$ the power $a^l$ has the
irreducible factorization
$$
X^l:=\{\langle\al,\,lX(\al)\rangle\cc\;
\al\in(X)_1\}\,.
$$

\begin{prop}[PP$\mathbf{0'}$]\label{prop_zerPrime}
Let $R$ be {\em UFD}, $a,b,c\in R$ and let $k\in\N$. If $a$
divides $bc^k$ and $a,c$ are coprime, then $a$ divides $b$.
\end{prop}
\duk
It is not hard to check that the proposition holds if $abc=0_R$. We  assume that $a,b,c\in R^*$ and denote by $X$, $Y$ and $Z$ their 
respective irreducible factorizations. Since $(X)_1\cap(Z)_1=\emptyset$, also 
$(X)_1\cap(Z^k)_1=\emptyset$. From the previous description of divisibility 
on the 
level of irreducible factorizations it then follows that since $X\,|\,YZ^k$, 
we in fact have $(X)_1\sus(Y)_1\setminus(Z)_1$ and $X\,|\,Y$. 
\kduk

\begin{prop}[PP$\mathbf{1'}$]\label{prop_PP3}
Let $R$ be {\em UFD} and let $k,l\in\N$ with $k,l\ge2$. If 
$a_i,b\in R$ for $i\in[k]$ are  
$k+1$ elements such that the $a_i$ are pairwise coprime and if
$$
a_1a_2\ds a_k\sim b^l\,,
$$
then there exist $k$ pairwise coprime elements $b_i\in R$ such that for every $i\in[k]$, 
$$
a_i\sim(b_i)^l\,.
$$
\end{prop}
\duk
Let $R$, $k$, $l$, $a_i$, and $b$ be as stated. If one of the $a_i$ 
is $0_R$, then every $a_j$ with $j\ne i$ is a~unit. The 
proposition then holds because $0_R=(0_R)^l$ and for every $a\in 
R^{\times}$ we have $a\sim(1_R)^l$. So we may omit every $a_i\in R^{\times}$ and may assume that 
$a_i,b\in R^*\setminus R^{\times}$. Let $X_i$ ($\ne\emptyset$) be the irreducible factorization of $a_i$, 
and $Y$ ($\ne\emptyset$) be that of $b$. Since
$$
X_1X_2\ds X_k=X_1\cup X_2\cup\ds\cup  X_k=Y^l\,,
$$
we have $X_i\sus Y^l$ for every $i\in[k]$. Thus $l$ divides $m$ for 
every $\langle\al,m\rangle\in X_i$, and it follows that $a_i\sim (b_i)^l$ for 
some $b_i\in R^*\setminus R^{\mathrm{ir}}$. Since the $a_i$ are pairwise coprime, so are the $b_i$.
\kduk

We use the following notation. If $R$ is a~domain and  
$A,B\sus R$ are two equinumerous finite sets, then $A\sim B$ means that there 
is a~bijection $f\cc A\to B$ such that $a\sim f(a)$ for every $a\in A$.

\begin{prop}[PP$\mathbf{2'}$]\label{prop_PP4}
Let $R$ be {\em UFD}, $p\in R^{\mathrm{ir}}$ and let $k\in\N$ with $k\ge2$. If
$a,b,c$ in $R$ are such that
$$
\mathrm{gcd}(a,\,b)=p\wedge ab\sim c^k\,,
$$
then there exist coprime elements $d,e\in R$ such that
$$
\{a,\,b\}\sim\{pd^k,\,p^{k-1}e^k\}\,.
$$
\end{prop}
\duk
Let $R$, $p$, $k$, $a$, $b$ and $c$ be as stated. If one of  
$a$ and $b$ is $0_R$, then the other 
is associated with $p$. We are done because $0_R=p^{k-1}(0_R)^k$ and 
$p=p(1_R)^k$. We may therefore assume that 
$a,b,c\in R^*$. Let $X$, $Y$ and $Z$ be the irreducible factorizations of 
$a$, $b$ and $c$, respectively. We may assume that $\langle[p]_{\sim},1\rangle\in X$, 
$\langle[p]_{\sim},m\rangle\in Y$ with $m\in\N$ and that $(X)_1\cap(Y)_1=\{[p]_{\sim}\}$. 
Since
$$
XY=\big(X\setminus\{\langle[p]_{\sim},\,1\rangle\}\big)
\cup\big(Y\setminus\{\langle[p]_{\sim},\,m\rangle\}\big)
\cup\{\langle[p]_{\sim},\,
m+1\rangle\}=Z^k\,,
$$
we get that $k$ divides both $X(\al)$ and $Y(\al)$ for every $\al\in((X)_1\cup 
(Y)_1)\setminus\{[p]_{\sim}\}$, and that $k$ divides $m+1$. Thus $m=k-1+km_0$ 
for some $m_0\in\N_0$. It follows that $a\sim pd^k$ and $b\sim p^{k-1}e^k$ 
for coprime $d,e\in R^*$.
\kduk

\noindent
Assuming that $\Z$ is UFD, we leave the deduction of PP0 from PP$0'$, PP1 from PP$1'$ and PP2 from PP$2'$ as easy exercises for the interested reader.

\section[${}^c$Klazar's resolution of $x^2-y^3=1$]{Klazar's resolution of $x^2-y^3=1$}\label{sec_klaz}

We follow the article \cite{klaz}. 
In the next result, we obtain well known formulas for Pythagorean triples. These formulas reverse 
the polynomial identity
$$
\big(x^2-y^2\big)^2+(2xy)^2=
\big(x^2+y^2\big)^2\ \ (\text{in }\Z[x,\,y])\,.
$$

\begin{prop}\label{prop_PythTrip}
If $x,y,z\in\N_0$ are numbers such that
$$
x^2+y^2=z^2
$$
---\,they form a~Pythagorean 
triple\,---\,and are pairwise
coprime, then they express for some coprime numbers
$u,v\in\N_0$ as 
$$
z=u^2+v^2\,\text{ and }\,
\{x,\,y\}=\{u^2-v^2,\,2uv\}\,.
$$
\end{prop}
\duk
Let $x$, $y$ and $z$ be as stated. Reduction modulo $4$ shows that $z$ and
exactly one of $x$ and $y$, say $x$, is odd. Then, since $(z-x)/2$ and $(z+x)/2$ are coprime, from the equality
$$
\Big(\frac{y}{2}\Big)^2=\frac{z-x}
{2}\cdot\frac{z+x}{2}
$$ 
we
get by PP1 (Proposition~\ref{prop_PP1}) coprime numbers $u,v\in\N_0$ such that $(z-x)/2=v^2$, 
$(z+x)/2=u^2$ and $uv=y/2$. Hence
$x=u^2-v^2$, $y=2uv$ and $z=u^2+v^2$.
\kduk

\begin{cor}\label{cor_onab}
If $x,y\in\N_0$ are numbers such that
$$
2x^2-y^2=1\,,
$$ 
then there exist numbers $a,b\in\N_0$ such that $a^2-2b^2=1$ and
$$
x\in\{a^2+2b^2+2ab,\,a^2+2b^2-2ab\}=\{2a^2-1+2ab,\,2a^2-1-2ab\}\,.
$$
\end{cor}
\duk
Let $x,y$ be as stated. Then $y$ is odd, $y=2y_0+1$ with $y_0\in\N_0$. So 
$x^2=2y_0^2+2y_0+1=y_0^2+(y_0+1)^2$. By Proposition~\ref{prop_PythTrip} there exist $u,v\in\N_0$ such that 
$$
x=u^2+v^2\,\text{ and }\,
\{y_0,\,y_0+1\}=\{2uv,\,u^2-
v^2\}\,.
$$
Thus $1=u^2-v^2-2uv=(u-v)^2-2v^2$ or $1=2uv-u^2+v^2=(u+v)^2-2u^2$. We set $a=u-v$,
$b=v$, respectively $a=u+v$, 
$b=u$, and get $a,b\in\N_0$ as stated.
\kduk

We took the following result from Sierpinski's book \cite{sier}.

\begin{prop}\label{prop_x4_2y2}
The only solutions of the Diophantine equation
$$
x^4-2y^2=1
$$
are $\langle\pm1,0\rangle$.
\end{prop}
\duk
Let $x,y\in\Z$ be a~solution of the equation. Then $x$ is odd and $x^2=1+4k$ for some $k\in\N_0$. We get 
from the factorization
$$
(x^2-1)(x^2+1)=2y^2
$$ 
that $4k(2k+1)=y^2$. Since the numbers $4k$ and $2k+1$ are coprime, by PP1 (Proposition~\ref{prop_PP1}) we have $4k=a^2$, $a\in\N_0$. Thus $(x-a)(x+a)=1$ and $x=\pm1,y=0$.
\kduk

The following auxiliary result is of an independent interest.

\begin{thm}\label{thm_x4_3y2}
The Diophantine equation
$$
x^4-3y^2=1
$$
has only the trivial solutions $\langle\pm1,0\rangle$.  
\end{thm}
\duk
We need to solve 
$x_n=m^2$, where $n,m\in\N_0$ and $x_n$ are as in Proposition~\ref{prop_x2_3y2}. If $n=2n_0+1$ is odd then 
Proposition~\ref{prop_x2_3y2Q} gives 
$$
(m-y_{n_0}-y_{n_0+1})(m+y_{n_0}+y_{n_0+1})=1\,. 
$$
Thus $m=\pm1$ and $y_{n_0}+y_{n_0+1}=0$. This is impossible because always $y_{n_0}+y_{n_0+1}>0$.

Let $n=2n_0$ be even. Then by Proposition~\ref{prop_x2_3y2Q}, $2x_{n_0}^2-1=x_{2n_0}=x_n=m^2$ and 
$2x_{n_0}^2-m^2=1$. Reduction modulo $4$ gives that $x_{n_0}$ is odd. By Proposition~\ref{prop_x2_3y2Q},
$x_{n_0}=x_{2n_1}=2x_{n_1}^2-1$. Corollary~\ref{cor_onab} shows that 
there are numbers $a,b\in\N_0$ such that $a^2-2b^2=1$ and
$$
2x_{n_1}^2-1=x_{n_0}=2a^2-1\pm2ab\,.
$$
Hence $x_{n_1}^2=a(a\pm b)$. Since $a,b$ are coprime, so are $a,a\pm b$,
and by PP1 (Proposition~\ref{prop_PP1}) the number $a$ is a~square. By
Proposition~\ref{prop_x4_2y2} we have $a=1$. Thus $b=0$, 
$x_n=x_{n_0}=1$ and $x=\pm1,y=0$. 
\kduk

\begin{thm}[Euler, 1737]\label{thm_x2y3}
The only integral solutions of the Diophantine equation
$$
x^2-y^3=1
$$
are $\langle\pm3,2\rangle$, $\langle\pm1,0\rangle$ and $\langle0,-1\rangle$. 
\end{thm}
\duk
(Klazar) Let $x,y\in\Z$ be such that $x^2-y^3=1$. Then 
$$
x^2=(y+1)\cdot(y^2-y+1)=(y+1)\cdot((y+1)(y-2)+3)\,.
$$
Since $y^2-y+1\ge0$, we have $y+1\ge0$. Also,  
$\mathrm{gcd}(y+1,y^2-y+1)\in\{1,3\}$. If the gcd is~1 we 
have by PP1 (Proposition~\ref{prop_PP1}) that $y+1$ and $y^2-y+1$ are squares.
Thus $4y^2-4y+4=(2a)^2$ for some $a\in\N_0$. We get that 
$$
3=(2a-2y+1)(2a+2y-1)\,.
$$
Hence $\langle a,y\rangle=\langle\pm1,1\rangle$ or 
$\langle\pm1,0\rangle$. For $y=1$ the number $y+1$ is not a~square
and for $y=0$ we get the solution $\langle\pm1,0\rangle$.

Let $\mathrm{gcd}(y+1,y^2-y+1)=3$. By PP2 (Proposition~\ref{prop_PP2}) there
are numbers $a,b\in\N_0$ such that 
$$
y+1=3a^2\,\text{ and }\,y^2-y+1=3b^2\,. 
$$
Thus $3(2b)^2-(2y-1)^2=3$ and $2y-1=3Y$
for some $Y\in\Z$. With $X=2b$ ($\in\N_0$) we get
$$
X^2-3Y^2=1\,\text{ and }\,Y=2a^2-1\,.
$$
The triple $\langle X,Y,a\rangle=\langle2,-1,0\rangle$ 
solves this system. We get $y=-1$ and the solution 
$\langle0,-1\rangle$.

Thus we can assume that $Y\ge0$. We look for a~number $n\in\N_0$
such that $Y=y_n=2a^2-1$, where $a\in\N_0$ and $y_n$ is as in
Proposition~\ref{prop_x2_3y2}. Since $y_n$ is odd, so is $n$. By Proposition~\ref{prop_x2_3y2Q} we have 
$$
2a^2-1=y_n=y_{2m+1}=2x_my_{m+1}-1,\ m\in\N_0\,.
$$
Thus 
$$
a^2=x_my_{m+1}\,.
$$
We have 
$$
\mathrm{gcd}(x_m,\,y_{m+1})=\mathrm{gcd}(x_m,\,x_m+2y_m)
\in\{1,\,2\}\,. 
$$
If $\mathrm{gcd}(x_m,y_{m+1})=1$, by PP1 (Proposition~\ref{prop_PP1}) the number $x_m$ is a~square.
Theorem~\ref{thm_x4_3y2} gives $x_m=1$. Thus $m=0$, $n=1$, $y_n=Y=1$
and $y=2$. We get the solution $\langle\pm3,2\rangle$. 

Finally, let $\mathrm{gcd}(x_m,y_{m+1})=2$. By PP2 (Proposition~\ref{prop_PP2}), $y_{m+1}=2c^2$ for some $c\in\N_0$. By
Proposition~\ref{prop_x2_3y2Q},
$$
2c^2=y_{m+1}=y_{2k}=2x_ky_k,\ k\in\N\,.
$$
Hence $c^2=x_ky_k$. Since $x_k,y_k$ are coprime, by PP1 (Proposition~\ref{prop_PP1}) the number $x_k$ is a~square. Theorem~\ref{thm_x4_3y2} shows that it is not
possible because $k\ge1$. We do not get any more solutions of $x^2-y^3=1$, and the proof 
is complete.
\kduk

\noindent
This was the 1989 resolution of $x^2-y^3=1$ in \cite{klaz}. In 2003, a~similar resolution was obtained by Ch.~Notari in \cite{nota}. 

\section[${}^c$Modern completion of Euler's argument]{Modern completion of Euler's argument}\label{sec_descent}

We give a~modern version of  Euler's resolution of the equation $x^2-y^3=1$ in the domain of 
fractions. We follow (and complete) the write-up \cite{conr} of Conrad.
First we resolve an auxiliary Diophantine equation.

\begin{thm}\label{thm_pomocnaV}
The only triple $\langle x,y,z\rangle\in\N^3$ such that $(x,y)=1$, $(3,x)=1$ and
$$
x^4-3x^2y^2+3y^4=z^2\,,
$$
is $\langle1,1,1\rangle$. Said equivalently, the only coprime numbers $u,v\in\N$ such that $(3,u)=1$ and that all three numbers
$$
u,\ v\,\text{ and }\,u^2-3uv+3v^2
$$
are squares, are $u=v=1$.
\end{thm}
\duk
(Conrad) First, we show that in the latter problem, if $v=1$ then $u=1$.
Indeed, since $u^2-3uv+3v^2$ is a~square, there is an $a\in\Z$ such that 
$$
u^2-3u+3=a^2\,\text{ and }\,(2a)^2-
(2u-3)^2=3\,.
$$
Thus $a=\pm1$ and $u=1$ or $2$. Since $u=2$ is not a~square, $u=1$.

We therefore assume that $u,v\in\N$ and $w\in\Z$ are numbers such that $(u,v)=1$, $(3,u)=1$, $u$ and $v$ are squares, 
$$
u^2-3uv+3v^2=w^2\,,
$$ 
$v>1$ and that $v$ is the minimum possible (with respect to the stated 
properties). We obtain a~contradiction and show that no such triple $\langle u,v,w\rangle$ in fact 
exists. 

Thus 
$(3,w)=1$ and we may select the sign of $w$ so that $w\equiv-u\ (\mathrm{mod}\;3)$. Let 
$r:=\frac{w-u}{v}$ ($\in\Q$), so that $u+rv=w$. By the choice of $w$ we have 
$r\ne0$. Let $r=\frac{m}{n}$ be in lowest terms, so that $m\in\Z$, 
$n\in\N$ and $(m,n)=1$. We see that $m\,|\,(w-u)$, $n\,|\,v$ and 
$(3,m)=1$. Substituting $w=u+rv$ in the previous displayed equation we get
$$
u^2-3uv+3v^2=u^2+2urv+r^2v^2\,\text{ and }\,(3-r^2)v^2=(2r+3)uv\,.
$$
We see that $2r+3\ne0$. Dividing by $v^2(2r+3)$ and substituting $r=\frac{m}{n}$ we get
$$
\frac{u}{v}=\frac{3n^2-m^2}{n(2m+3n)}\,.
$$
We show that the fraction on the right-hand side is in lowest terms. From 
$(m,n)=1$ we get that $n$ is coprime with $3n^2-m^2$. To show that $3n^2-
m^2$ and $2m+3n$ are coprime, we argue by contradiction. Suppose that a~prime 
$p$ divides $3n^2-
m^2$ and $2m+3n$. Then $m^2\equiv3n^2$ and $2m\equiv-3n$
modulo $p$. Squaring the second congruence and comparing it with the 
first we get $4m^2\equiv3m^2$ and $12n^2\equiv9n^2$. Thus $p$ divides 
both $m^2$ and $3n^2$. Since $m$ and $n$ are coprime, $p=3$ and $p$ divides
$m$, which is a~contradiction.

Thus
$$
u=\ep(3n^2-m^2)\,\text{ and }\,v=\ep n(2m+3n),\ \ep\in\{-1,\,1\}\,.
$$
Modulo $3$, $u\equiv-\ep m^2\equiv-\ep$. Since $u$ is a~square, we get that $\ep=-1$. Hence 
$$
u=m^2-3n^2\,\text{ and }\,v=-n(2m+3n)\,.
$$
Since $u$ is a~square, we have
$$
m^2-3n^2=k^2,\ k\in\Z\,.
$$
Thus $(3,k)=1$ and we select the sign of $k$ so that $k\equiv-m$ modulo $3$.

Let $s=\frac{k-m}{n}$ ($\in\Q$), so that $m+sn=k$. It follows from the 
choice of $k$ that $s\ne0$. Let $s=\frac{u'}{v'}$ be in lowest terms, so that $u'\in\Z$, 
$v'\in\N$ and $(u',v')=1$. We see that $u'\,|\,(k-m)$, $v'\,|\,n$ and 
$(3,u')=1$. Substituting $k=m+sn$ in the previous displayed equation we get
$$
m^2-3n^2=m^2+2mns+s^2n^2\,\text{ and }\,2mns=-(3+s^2)n^2\,.
$$
We know that $sn\ne0$. Dividing by $sn^2$ and substituting $s=\frac{u'}{v'}$ we get
$$
\frac{2m}{n}=-\frac{3+s^2}{s}\,\text{ and }\,v=-n^2\Big(\frac{2m}{n}+3\Big)=
n^2\cdot\frac{(u')^2-3u'v'+3(v')^2}{u'v'}\,.
$$
Since $v$ is a~square, multiplying by $(u'v')^2$ we get that
$$
u'\cdot v'\cdot \big((u')^2-3u'v'+3(v')^2\big)
$$
is a~square. Since $(u',v')=1$ and $(3,u')=1$, the three displayed 
factors are pairwise coprime and by PP1 (Proposition~\ref{prop_PP1}) all 
three are squares; note that $v'>0$ 
and $(u')^2-3u'v'+3(v')^2>0$, hence also $u'>0$.

We show that $0<v'\le v$. We have 
$$
\frac{v}{n}=\frac{n\big((u')^2-3u'v'+3(v')^2\big)}{u'v'}\,.
$$
Now $n\,|\,v$, and $u'$ and $v'$ are coprime with $(u')^2-3u'v'+3(v')^2$. 
Thus $u'v'\,|\,n$ and $u'v'\,|\,v$. Hence $0<u'v'\le v$ and 
$0<v'\le v$. The minimality of $v$ implies two cases: 
either $v'=1$ or $v'=v$. We show that both lead to the contradiction that $v=1$. 

{\bf Case~1 when $v'=1$. }As we know from the beginning, $u'=1$. Thus 
$s=\frac{u'}{v'}=1$, $k=m+sn=m+n$, $m^2-3n^2=k^2=(m+n)^2$ and 
$m=-2n$. So $r=\frac{m}{n}=-2$ and $w=u+rv=u-2v$. We get
$$
u^2-3uv+3v^2=(u-2v)^2\,\text{ and }\,
uv=v^2\,.
$$
Hence $(u-v)v=0$ and $u=v$. Since $(u,v)=1$, we get the contradiction that $v=1$.

{\bf Case~2 when $v'=v$. }Thus $u'=1$.
We have $s=\frac{k-m}{n}=\frac{u'}{v'}=\frac{1}{n}$. Since $n\,|\,v$ and $v'\,|\,n$, we get $v=v'=n$ and $k=m+1$. Then $m^2-3n^2=k^2=m^2+2m+1$ and $2m+1=-3n^2$.
But from
$$
n=v=-n(2m+3n)
$$
we get that $2m+3n=-1$ and $2m+1=-3n$.
Thus $-3n^2=-3n$ and $n=v=1$, which is a~contradiction.   
\kduk

\noindent 
In the previous proof we followed \cite{conr}, but we 
corrected/completed the conclusion of the proof; \cite{conr} misses case~1. 

From the previous theorem we deduce the main result of this section

\begin{thm}[Euler, 1737]\label{thm_eulerDesc}
The only rational solutions $x,y\in\Q$ of the equation
$$
x^2-y^3=1
$$
are the pairs $\langle\pm3,2\rangle$, $\langle\pm1,0\rangle$ and $\langle0,-1\rangle$.
\end{thm}
\duk
(Conrad) Let $x,y\in\Q$ be such that $x^2-y^3=1$. Since $x^2,y^2-y+1=(y-\frac{1}{2})^2+\frac{3}{4}\ge0$, from 
$$
x^2=y^3+1=
(y+1)(y^2-y+1)
$$ 
we deduce that $y\ge-1$. As $y=-1$ yields the solution 
$\langle0,-1\rangle$, we assume from now on that $y>-1$.  

Let $y=\frac{a}{b}$ with coprime $a\in\Z$ and $b\in\N$.  From $y>-1$ we get $a+b>0$. Since 
$$
x^2=y^3+1=\frac{b(a^3+b^3)}{b^4}\,,
$$
$b(a^3+b^3)$ is a~square. With $c=a+b$ ($\in\N$) we write
$$
b(a^3+b^3)=b(a+b)(a^2-ab+b^2)=
b\cdot c\cdot(c^2-3bc+3b^2)\,.
$$

{\bf Case~1 when $(3,c)=1$. }The three displayed factors are pairwise coprime 
and positive. By PP1 
(Proposition~\ref{prop_PP1})  they are squares. Using 
Theorem~\ref{thm_pomocnaV} we get that $b=c=1$ and $a=0$. Thus $y=\frac{a}{b}=0$ and we get the solution $\langle\pm1,0\rangle$.

{\bf Case~2 when $3\,|\,c$. }Let $c=3d$. From $(b,c)=1$ we get $(3,b)=1$
and $(b,d)=1$. Then 
$$
bc(c^2-3bc+3b^2)=9bd(b^2-3bd+3d^2)
$$
is a~square.  Thus
$$
b\cdot d\cdot(b^2-3bd+3d^2)
$$
is a~square. As before we deduce by means of PP1 (Proposition~\ref{prop_PP1}) and 
Theorem~\ref{thm_pomocnaV} that $b=d=1$. Thus $c=3d=3$, $a=c-b=2$ and
$y=\frac{a}{b}=2$. We get the solution $\langle\pm3,2\rangle$.
\kduk

\section{Wakulicz's resolution of $x^3+y^3=2z^3$}\label{sec_wakulicz}

Antoni Wakulicz (1902--1988) (\cite{wakuliczWiki})

\chapter[V.~Lebesgue's theorem:  $x^m-y^2=1$]{V.~Lebesgue's theorem}

In the second chapter, we tackle the second elementary case of Catalan's conjecture, that of equations 
$$
x^m-y^2=1
$$
with odd $m\ge3$ (for even $m$, it is a~trivial problem). Section~\ref{sec_pAdic} contains results on $p$-adic order that
we use later in Section~\ref{sec_lebe}, where we present the theorem of V.~Lebesgue (1791--1875) that the only solution of 
the equation is $\langle1,0\rangle$. For a~biography of V.~Lebesgue, see 
\cite{Vlebe}. In Section~\ref{sec_eucUFD} we prove in 
Theorem~\ref{thm_euclIsUFD} that every Euclidean domain is UFD. In 
particular, $\Z[i]$ 
is UFD, which is needed for 
Lebesgue's proof. This theorem plays a more important role in the resolution of Catalan's 
conjecture than one might think. We show by means of it that the class 
numbers of the cyclotomic fields $\Q(\zeta_3)$ and $\Q(\zeta_5)$ are
$h_3=h_5=1$, which is a~key result in Chapter~\ref{chap_M4}. In Section~\ref{sec_lebe} we follow \cite[Chapter~2]{scho}. 

\section[${}^c$Properties of $p$-adic order]{Properties of $p$-adic order}\label{sec_pAdic}

For a prime $p$ and $\al\in\Q$, the $p$-adic order 
$\mathrm{ord}_p(\al)$ of $\al$ is $+\infty$ if $\al=0$. For $\al\ne0$, it is the unique $k\in\Z$
such that $\al=p^k\be$ for some $\be\in\Q$ with numerator and 
denominator coprime to $p$.

\begin{prop}\label{prop_vlPord}
For every $\al,\be\in\Q$, the following holds.
\begin{enumerate}
\item $\mathrm{ord}_p(\al\be)=\mathrm{ord}_p(\al)+\mathrm{ord}_p(\be)$. Here for $a,b\in\Z\cup\{+\infty\}$ we set $a+b:=+\infty$ if $a$ or $b$ is $+\infty$.
\item $\mathrm{ord}_p(\al+\be)\ge\min(\mathrm{ord}_p(\al),\mathrm{ord}_p(\be))$, with equality if $\mathrm{ord}_p(\al)\ne\mathrm{ord}_p(\be)$. Here $+\infty>m$ for every $m\in\Z$. 
\end{enumerate}
\end{prop}
\duk
1. If $\al$ or $\be$ is $0$ then the equality holds, and we may assume that $\al,\be\ne0$. 
Let $\al=p^k\al_0$ and $\be=p^l\be_0$ in $\Q^*$ be two 
arbitrary fractions, where $k=\mathrm{ord}_p(\al)$,
$l=\mathrm{ord}_p(\be)$ and 
the fractions $\al_0$ and $\be_0$
have numerators and denominators coprime to $p$. Then
$$
\al\be=p^{k+l}\al_0\be_0=:p^{k+l}
\gamma
$$
where $\gamma\in\Q^*$
has numerator and denominator coprime to $p$. Thus 
$$
\mathrm{ord}_p(\al\be)=k+l=
\mathrm{ord}_p(\al)+\mathrm{ord}_p(\be)\,.
$$

2. Again, if $\al$ or $\be$ is $0$ then the claim holds.
We take $\al,\be\in\Q^*$ as in item~1 and assume without loss of generality that $k\le l$. Then
$$
k=\min(\mathrm{ord}_p(\al),\,\mathrm{ord}_p(\be))
$$
and 
$$
\al+\be=p^k(\al_0+p^{l-k}\be_0)=:
p^k\gamma\,.
$$
If $k<l$ then $\gamma$ can be written as a~fraction 
with numerator and denominator coprime to $p$, so that 
$\mathrm{ord}_p(\al+\be)=k$. If $k=l$ then $\gamma$ can 
be written with denominator coprime to $p$, and it follows that 
$\mathrm{ord}_p(\al+\be)\ge k$.
\kduk

\begin{cor}\label{cor_pAdicL} 
If $p$ is a~prime and $\al_1$, 
$\dots$, $\al_n$ are $n\ge2$ fractions such that
$$
\mathrm{ord}_p(\al_n)<\mathrm{ord}_p(\al_i)\,\text{ for 
every $i\in[n-1]$}\,,
$$
then $\sum_{j=1}^n\al_j\ne0$. 
\end{cor}
\duk
Let $k:=\min(\{\mathrm{ord}_p(\al_i)\colon\;i\in[n-1]\})$ ($\in\Z\cup\{+\infty\}$) and let 
$\al:=\sum_{j=1}^{n-1}\al_j$. Applying repeatedly item~2 of 
Proposition~1 we get that $\mathrm{ord}_p(\al)\ge k$. Since $\mathrm{ord}_p(\al_n)<k$, we have $\mathrm{ord}_p(\al_n)<+\infty$ and
$${\textstyle
\mathrm{ord}_p\big(\sum_{j=1}^n \al_j\big)=\mathrm{ord}_p(\al+\al_n)=
\mathrm{ord}_p(\al_n)<+\infty\,.
}
$$
Thus $\sum_{j=1}^n\al_j\ne0$.
\kduk

\section[${}^c$V.~Lebesgue's resolution of $x^m-y^2=1$]{V.~Lebesgue's resolution of $x^m-y^2=1$}\label{sec_lebe}

We work in the domain
$$
\Z[i]=\langle
\Z[i],\,0,\,1,\,+,\,\cdot,\,\rangle\,,
$$
where the base set $\Z[i]=
\{a+bi\cc\;a,b\in\Z\}$ is the set of complex numbers with both real 
and imaginary part integral.

\begin{thm}[V.~Lebesgue, 1850]\label{thm_LebeThm}    
Let $m\ge3$ be an odd integer. The Diophantine equation 
$$
x^m-y^2=1
$$
has only the solution $\langle 1,0\rangle$.
\end{thm}
\duk
(V.~Lebesgue) Let $m$ be as stated and $a,b\in\Z$ with $b\ne0$ be such that 
$a^m-b^2=1$. We derive a~contradiction. If $b$ is odd then $a^m\equiv2$ modulo $4$, which is 
impossible. Thus $b$ is even, $b\ne0$ and $a$ is odd. In $\Z[i]$, we
consider the factorization
$$
a^m=(1+bi)(1-bi)\,.
$$
The elements $1+bi$ and $1-bi$ are coprime: if $\al\in\Z[i]$ divides both $1+bi$ and $1-bi$, then $n=\al\overline{\al}$ ($\in\N$) divides, in $\Z$, the 
number $2\cdot\overline{2}=4$ and the odd number $(1+bi)(1-bi)=a^m$. Thus 
$n=1$, $\al$ is a~unit, and $1+bi$ 
and $1-bi$ are coprime. Since $\Z[i]$ is UFD (Corollary~\ref{cor_ziUFD} in the next section), by PP$1'$ 
(Proposition~\ref{prop_PP3})
there exist an element $\al\in\Z[i]$, units
$\epsilon,\epsilon'\in\Z[i]^{\times}$ 
and numbers $u,v\in\Z$ such that
$$
1+bi=\epsilon\al^m=(\epsilon'\al)^m=(u+vi)^m\,\text{ and }\,
1-bi=\overline{\epsilon}(\overline{\al})^m=(\overline{\epsilon'\al})^m=(u-vi)^m\,,
$$
because every unit
in $\Z[i]^{\times}=\{\pm1,\pm i\}$ (Proposition~\ref{prop_ziUnits} in the next section) is an $m$-th power. Since $m$ is odd, we have
$$
2=(u+vi)^m+(u-vi)^m=2u\cdot\be,\ \be\in\Z[i]\,,
$$
and deduce that $u=\pm1$. We exclude the possibility $u=-1$. Since 
$(1+v^2)^m=(u^2+v^2)^m=1+b^2$ is odd, the number $v$ is even. From
$${\textstyle
1+bi=(u+vi)^m=\sum_{j=0}^m
\binom{m}{j}u^{m-j}(vi)^j\equiv u^m+
mu^{m-1}vi\ (\mathrm{mod}\;4)
}
$$
(congruence in $\Z[i]$) we deduce 
that $u^m\equiv1$ modulo $4$ (congruence in $\Z$), which excludes $u=-1$.

Thus $u+vi=1+vi$ with even and nonzero $v$ (since $b\ne0$). 
Comparing the real parts in $1+bi=(1+vi)^m$ we get an identity in 
$\Z$,
$$
{\textstyle
1=\sum_{j=0}^{(m-1)/2}(-1)^j\binom{m}{2j}v^{2j},\,\text{ or }\,
-\binom{m}{2}v^2+\sum_{j=2}^{(m-1)/2}(-1)^j\binom{m}{2j}v^{2j}=0\,.
}
$$
For $m=3$ the last sum is empty (zero) and the equality is impossible 
as $v\ne0$. For odd $m\ge5$ we show that the equality does not hold by 
means of Corollary~\ref{cor_pAdicL} and prime $p=2$. We set $A=\binom{m}{2}v^2$ and 
$B_j=\binom{m}{2j}\cdot v^{2j}$ for $j=2,3,\ds,\frac{m-1}{2}$, and show
that $\mathrm{ord}_2(A)<\mathrm{ord}_2(B_j)$ for every $j$. 
Indeed,
$$
{\textstyle
B_j=A\cdot\frac{1}{j(2j-1)}\binom{m-2}{2j-2}v^{2j-2}=:A\cdot C_j
}
$$
and $\mathrm{ord}_2(C_j)\ge 2j-2-\lfloor\log_2(j)\rfloor>0$, so that
by the additivity of $\mathrm{ord}_2(\cdot)$ (item~1 of 
Proposition~\ref{prop_vlPord}) we have $\mathrm{ord}_2(A)=
\mathrm{ord}_2(B_j)-
\mathrm{ord}_2(C_j)<\mathrm{ord}_2(B_j)$. We get 
a~contradiction
\kduk

\section[${}^c$Euclidean domains are UFD]{Euclidean domains are UFD}\label{sec_eucUFD}

We begin by reviewing Euclidean domains. A~{\em well ordering} 
$\langle W,\prec\rangle$ is a~linear
order $\prec$ on a~set $W$ such that every nonempty set $V\sus W$ has
the {\em minimum element $m\in V$}, an element such that $m\preceq x$ for 
every $x\in V$. Minima are unique.

\begin{defi}[Euclidean domain]\label{def_euclDom}
A~domain 
$$
R=\langle R,\,0_R,\,1_R,\,+,\,\cdot\rangle
$$ 
is called Euclidean if there is a~well ordering $\langle 
W,\prec\rangle$ and a~function 
$f\cc R^*\to W$ such that
$$
\forall\,a,\,b\in R,\,b\ne 0_R\,
\exists\,c,\,d\in R\,\big(
a=b\cdot c+d\wedge(d=0_R\vee f(d)\prec f(b))\big)\,.
$$
\end{defi}
Note that the last disjunction is an exclusive or. For
instance, the domain of integers $\Z$ is Euclidean: 
$\langle W,\prec\rangle$ is
$\langle \N,<\rangle$ and $f(n)=|n|$. One can often prove that a~domain is 
UFD by using the next classical theorem.

\begin{thm}\label{thm_euclIsUFD}
Every Euclidean domain is {\em UFD}.  \end{thm}
\duk
Let $R$ be a~Euclidean domain with a~well ordering 
$\langle W,\prec\rangle$ and a~map $f\cc R^*\to W$ as in Definition~\ref{def_euclDom}. 

{\em Existence of irreducible factorizations. }We show that every 
$x\in R^*\setminus R^{\times}$ is a~product of irreducibles. Suppose for 
the contrary
that the set 
$$
A\sus R^*\setminus R^{\times}
$$ 
of elements that are not
products of irreducibles is nonempty. Let $a\in R^*$ be
such that $a$ has a~divisor $b\in A$,  and that the value $f(a)$ is $\prec$-minimum among
all such values in $W$. Thus
$a=bc$ where $b\in A$ and $c\in R^*$. Since $b$ is not irreducible, $b=de$
with $d,e\in R^*\setminus R^{\times}$. But
$b\in A$ and hence $d$ or $e$ is in $A$. We assume that $d\in A$, the case
with $e\in A$ is similar. Thus
$a=d(ec)$ where $d\in A$ and $ec\in R^*\setminus R^{\times}$ (because
$e\in R^*\setminus R^{\times}$). This means by Proposition~\ref{prop_mutuDiv} that 
$a$ does not divide $d$, and if we divide $d$ by $a$ with a~remainder we get
$$
d=ag+h\,\text{ where $g\in R$, $h\in R^*$ and $f(h)\prec f(a)$}\,.
$$
Since $d$ divides $a$, it divides $h$ too. So $d\in A$ and divides $h$, and $f(h)\prec f(a)$. This
contradicts the choice of $a$.

{\em Bachet's identity. }We prove that if $a,b\in R$ are coprime, then there
exist $c,d\in R$ such that $ca+db=1_R$. We consider the set
$$
I=\{ca+db:\;c,\,d\in R\}\ (\sus R)\,,
$$
which is the ideal in $R$ generated by the
elements $a,b$. Let $e\in I\setminus
\{0_R\}$ have $\prec$-minimum value $f(e)$. Clearly, $I\ne\emptyset$. Note that $I\ne\{0_R\}$ because we do not 
have $a=b=0_R$, the element $0_R$ is not a~unit and therefore $0_R,0_R$ are not coprime. We show that
$e$ divides every $x\in I$. Indeed, we express any $x\in I$ as $x=ec+d$ where $c,d\in R$ and $d=0_R$ or
$f(d)\prec f(e)$. Due to $d=x-ec\in I$ we have $d=0_R$. Thus $e$ divides every element of $I$ and since 
$a,b\in I$ and are coprime, $e\in R^{\times}$. It follows that $1_R\in I$, there exist $c,d\in R$ such that 
$1_R=ca+db$.

{\em Primes are prime divisors. }We show that if $a,b,c\in R$, $a\,|\,bc$ and $a\in R^{\mathrm{ir}}$, then $a\,|\,b$ or
$a\,|\,c$. Suppose that $a,b,c\in R$ and that $a$ is irreducible, divides
$bc$ but does not divide $b$. Then $a,b$ are coprime and by the second step there are $d,e\in R$ such that
$$
da+eb=1_R\,.
$$
We multiply it by $c$ and get $dac+ebc=c$. Hence $a$ divides $c$.

{\em Finally,} we prove that in $R$ every element 
$a\in R^*$ has 
a~unique irreducible 
factorization $X$ ($\sus\mathbb{I}\times\N$). It means to prove that if 
$$
a_1a_2\ds a_k\sim b_1b_2\ds b_l\,,
$$
where $a_i,b_i\in R^{\mathrm{ir}}$ and $k,l\in\N$, then always $k=l$ and there always exists 
a~bijection $f\cc[k]\to[l]$ such that  $$
a_i\sim b_{f(i)},\ i\in[k]\,.
$$ 
Suppose that the above displayed equality is a~shortest counterexample.
Then $k,l\ge2$ and, by the previous step, $a_k\sim b_m$ for some 
$m\in[l]$. Canceling $a_k$ and $b_m$ we get that
$$
a_1a_2\ds a_{k-1}\sim  b_1b_2\ds b_{m-1}b_{m+1}\ds b_l\,.
$$
Now we have $k-1=l-1$ and a~map $g\cc[k-1]\to[l]\setminus\{m\}$ with 
the above property. But then $k=l$ and we easily extend $g$ to $f\cc[k]\to[l]$ so that $f$
has the above property. We get a~contradiction and deduce that no 
counterexample actually exists.
\kduk

\begin{cor}\label{cor_ziUFD}
The domain $\Z[i]$ is Euclidean and hence {\em UFD}.    
\end{cor}
\duk
Recall that for $z=u+vi$ ($\in\C$), the norm (absolute value) of $z$ is 
$|z|=\sqrt{z\cdot\overline{z}}=
\sqrt{u^2+v^2}$. We have $|z\cdot z'|=|z|\cdot|z'|$ and 
$|z+z'|\le|z|+|z'|$. We take the well ordering $\langle\N,<\rangle$, i.e. 
the usual linear order on natural 
numbers, and the map
$$
f\cc\Z[i]^*\to\N,\ f(z)=|z|^2\,.
$$
Let $z,z'\in\Z[i]$ with $z'\ne0$. We define $\al,\be\in\Q$ by
$$
\frac{z}{z'}=\al+\be i\,.
$$
Let $a,b\in\Z$ be such that $|a-\al|\le\frac{1}{2}$ and $|b-\be|\le\frac{1}{2}$. Let $w:=a+bi$ and 
$w'=z-z'\cdot w$. Then $z=z'\cdot w+w'$, of course, and 
$$
f(w')=|w'|^2=|z'|^2\cdot\Big|
\frac{z}{z'}-w\Big|^2\le f(z')
\Big(\frac{1}{4}+\frac{1}{4}\Big)
<f(z')\,.
$$
\kduk

\begin{prop}\label{prop_ziUnits}
$\Z[i]^{\times}=\{-1,1,-i,i\}$.    
\end{prop}
\duk
Let $f(z)=|z|^2\cc\Z[i]^*\to\N$ be 
the map in the previous proof. If $z=a+bi,z'\in\Z[i]$ are such that $zz'=1$, then
$$
1=f(1)=f(z)f(z')=(a^2+b^2)f(z')
$$
and $a^2+b^2=1$. Thus $a=\pm1, b=0$ or $a=0,b=\pm1$, and we get the stated units.
\kduk

\chapter[Chao Ko's theorem: $x^2-y^q=1$ with $q\ge5$]{Chao Ko's theorem}\label{chap_chaoKo}

The third chapter is devoted to the last elementary case of Catalan's 
conjecture, that of the equation 
$$
x^2-y^q=1 
$$
in which $q\ge5$ is a~prime 
number. In 1965, the Chinese mathematician Chao Ko (1910--2002) 
(\cite{chaoKoWiki}) proved in \cite{ko} that for every $q$ the only solutions are 
$\langle\pm1,0\rangle$ and $\langle0,-1\rangle$. Another 
well-known mathematical result due to Chao Ko is the Erd\H{o}s--Ko--Rado theorem in extremal combinatorics
\cite{EKR,EKRwiki}. 

Section~\ref{sec_2lem} contains two auxiliary lemmas. In 
Section~\ref{sec_chein} we present
a~resolution of the equation due to E.\,Z. Chein \cite{chein}. In this chapter we follow 
\cite[Chapter~3]{scho}.

\section[${}^c$Two lemmas]{Two lemmas}\label{sec_2lem}

\begin{lemma}\label{lem_oGCD}
Let $q$ be a~prime number and $a,b\in\Z$, $a\ne b$, be coprime numbers. Then   
$$
{\textstyle
d:=\mathrm{gcd}\big(\frac{a^q-b^q}{a-b},\,a-b\big)=\mathrm{gcd}\big(\sum_{i=0}^{q-1}a^ib^{q-1-i},\,a-b\big)
}
$$
divides $q$.
\end{lemma}
\duk
By the binomial theorem, $\frac{a^q-b^q}{a-b}=\frac{(a-b+b)^q-b^q}{a-b}$ equals
$$
{\textstyle
\sum_{i=1}^q\binom{q}{i}(a-b)^{i-1}b^{q-i}=qb^{q-1}+
\sum_{i=2}^q\binom{q}{i}(a-b)^{i-1}b^{q-i}\,.
}
$$
Thus $d\,|\,qb^{q-1}$. Since $a,b$ are coprime, so are $b,a-b$ and $(d,b)=1$. 
By PP0 (Proposition~\ref{prop_aDelBc}) the number $d$ divides $q$. 
\kduk

\begin{lemma}\label{lem_druheLem}
Let $q\in\N$ with $q\ge3$ be odd and $x,y$ be nonzero integers such that $x^2-y^q=1$. Then $y$ is even and
replacing $x$ with $-x$ if necessary, we have expressions
$$
x-1=2^{q-1}a^q\,\text{ and }\,x+1=2b^q\,,
$$
where $a,b\in\Z$ are coprime and $b$
is odd. 
\end{lemma}
\duk
In the factorization $(x-1)(x+1)=y^q$ the two factors are coprime or their gcd is $2$. In the former case we get by PP1 (Proposition~\ref{prop_PP1}) two 
$q$-th powers differing by $2$. The only such powers are $-1$ and $1$, and 
$x=0$. Since this is impossible, we see that $\mathrm{gcd}(x-1,x+1)=2$ and $y$ is even. Hence $x$ is odd. Changing 
the sign of $x$ we may assume that $\frac{x+1}{2}$ is odd. Applying PP2 
(Proposition~\ref{prop_PP2}) we get the stated expressions.
\kduk

\section[${}^c$Chein's resolution of $x^2-y^q=1$]{Chein's resolution of $x^2-y^q=1$}\label{sec_chein}

\begin{prop}\label{prop_prvTvrz}
Let $q\ge3$ be prime and $x,y$ be nonzero integers such that $x^2-y^q=1$. Then $q$ divides $x$.    
\end{prop}
\duk
We may assume that $y\in\N$. 
We assume that
$\neg(q\,|\,x)$ and obtain a~contradiction. By
Lemma~\ref{lem_oGCD}, in the factorization
$${\textstyle
x^2=(y+1)\cdot\frac{y^q-(-1)}{y-(-1)}
}
$$
the gcd of the two factors divides $q$. By the assumption on $x$ they are therefore coprime. 
By PP1 (Proposition~\ref{prop_PP1}) we have $y+1=u^2$ for $u\in\N$. Since $y$ is even, $u$ is odd. From the equalities
$$
x^2-y\cdot\big(y^{(q-1)/2}\big)^2=1\,
\text{ and }\,u^2-y\cdot1^2=1
$$
we see that $\langle x,y^{(q-1)/2}\rangle$ and $\langle u,1\rangle$ are two solutions of Pell equation
$$
X^2-yY^2=1\,. 
$$
(Since $y=u^2-1\ge3$ and is not a~square, it is really
a~Pell equation.) It is clear that $\langle u,1\rangle\in\N^2$ is the minimal solution. By
Proposition~\ref{prop_PellEq} there is an $m\in\N$ such that in the
domain $\Z[\sqrt{y}]$ we have the equality 
$$
x+y^{(q-1)/2}\sqrt{y}=\big(u+\sqrt{y}\big)^m\,.
$$
Thus in $\Z[\sqrt{y}]$ we have the congruence 
$$
x\equiv u^m+mu^{m-
1}\sqrt{y}\ (\mathrm{mod}\;y)\,.
$$
Hence, in $\Z$, the number $y$
divides $mu^{m-1}$. But $y$ is even and $u$ is odd, so $m$ is even. In $\Z[\sqrt{y}]$ we therefore have equality 
$$
x+y^{(q-1)/2}\sqrt{y}=\big(u^2+y+
2u\sqrt{y}\big)^{m/2}\,,
$$
from which we get the congruence $x+y^{(q-1)/2}\sqrt{y}\equiv y^{m/2}\ (\mathrm{mod}\;u)$. Thus in $\Z$ the
number $u$ divides $y^{(q-1)/2}$. But $y+1=u^2$, so $y,u$ are coprime and 
$u=\pm 1$. Hence $y=0$, in contradiction with the assumption.
\kduk

\begin{prop}\label{prop_x2yq_2}
Let $q\ge3$ be prime and $x,y$ be nonzero integers such that $x^2-y^q=1$. Then $x\equiv\pm3$ $(\mathrm{mod}\;q)$.   
\end{prop}
\duk
Let $q$, $x$ and $y$ be as stated. Changing the sign of
$x$ if necessary, by Lemma~\ref{lem_druheLem} we have coprime $a,b\in\Z$ with odd $b$ such that $x-1=2^{q-1}a^q$ and $x+1=2b^q$. Hence
$$
b^{2q}-(2a)^q=\Big(\frac{x+1}{2}\Big)^2-2(x-1)
=\Big(\frac{x-3}{2}\Big)^2
$$
and
$$
(b^2-2a)\cdot\Big(\frac{b^{2q}-(2a)^q}{b^2-2a}\Big)=
\Big(\frac{x-3}{2}\Big)^2\,.
$$
The numbers $2a,b^2$ are coprime and by Lemma~\ref{lem_oGCD}
the gcd of the last two factors divides $q$.

If it is $q$ then $x\equiv3$ $(\mathrm{mod}\;q)$ and for the original $x$, before
the possible change of sign, we have $x\equiv\pm3$ $(\mathrm{mod}\;q)$. We assume for the contrary that
the gcd is $1$. Hence we assume that $b^2-2a$ and $\frac{b^{2q}-(2a)^q}{b^2-2a}$
are coprime numbers. From $b^{2q}-(2a)^q\ge0$ (it is a~square) we get $b^2-2a\ge0$ ($f(X)=X^q$ is an increasing function) and by
PP1 (Proposition~\ref{prop_PP1}) there is a~$c\in\N$ such that $b^2-2a=c^2$. Since $y\ne0$, also $a\ne0$
and $c^2\ne b^2$. The nearest squares to $b^2$ different from it are $(b\pm1)^2$. Thus 
$2|a|=|b^2-c^2|\ge2|b|-1$ and hence $|a|\ge|b|$. On the other hand,
$$
|a|^q=\frac{|x-1|}{2^{q-1}}\le\frac{|x-1|}{16}<
\frac{|x+1|}{2}=|b|^q\,.
$$
For $x\in\Z$ the crucial strict inequality $|x-1|<8|x+1|$ does not hold only for $x=-1$. This value
of $x$ is excluded by the fact that $y\ne0$. Hence also $|a|<|b|$ and we have a~contradiction.
\kduk

\begin{thm}[Chao Ko, 1965]\label{thm_x2yq}
Let $q\ge5$ be a~prime number. The Diophantine equation
$$
x^2-y^q=1
$$
has only the solutions $\langle\pm1,0\rangle$ and 
$\langle0,-1\rangle$.
\end{thm}
\duk
(Chein) Suppose that $q$ is as stated and $x,y\in\Z^*$ satisfy $x^2-y^q=1$. By Proposition~\ref{prop_prvTvrz}
we have $x\equiv0\ (\mathrm{mod}\;q)$. Also $x\equiv\pm3\ 
(\mathrm{mod}\;q)$ by Proposition~\ref{prop_x2yq_2}. For $q>3$ these congruences
are contradictory.
\kduk

\chapter[Two relations of Cassels: $p\,|\,y$ and $q\,|\,x$]{Two relations of Cassels}

In Sections~\ref{sec_qdivx}  and \ref{sec_pdivy} we deduce two divisibility relations for hypothetical nonzero solutions $x,y\in\Z^*$ of the equation
$$
x^p-y^q=1\,,
$$
where $p>q>2$ are primes: $q$ divides $x$ and $p$ divides $y$. 
Since we eventually show that no such numbers $x$ and $y$ exist, these are 
properties of non-existing objects. These relations, important in 
Mih\u{a}ilescu's proof, were obtained in \cite{cass1} and \cite{cass2}, respectively, by the 
British mathematician John~W.~S.~Cassels (1922--2015) 
(\cite{casselsWiki}). Section~\ref{sec_fiveLem} contains 
auxiliary results. In 
Section~\ref{sec_corrCass} we obtain corollaries of the two relations. 
In this chapter we 
follow \cite[Chapter~6]{scho}. 

\section{Five lemmas}\label{sec_fiveLem}

\begin{lemma}\label{lem_incrDecr}
For $u\in\R$ the following hold. 
\begin{enumerate}
\item If $u\ge1$ then the function 
$$
f(x)=\big(u^x+1\big)^{1/x}\cc(0,\,
+\infty)\to(0,\,+\infty)
$$ 
decreases. 
\item If $u>1$ then the function 
$$
f(x)=\big(u^x-1\big)^{1/x}\cc(0,\,+\infty)\to(0,\,
+\infty)
$$ 
increases. 
\end{enumerate}
\end{lemma}
\duk
1. We have 
$$
f'(x)=
f(x)\Big(\frac{u^x\log u}{x(u^x+1)}-\frac{\log(u^x+1)}{x^2}\Big)
$$ 
and $(\cdots)<0$ because $xu^x\log u=u^x\log(u^x)<(u^x+1)\log(u^x+1)$. 

2. Similarly, 
$$
f'(x)=
f(x)\Big(\frac{u^x\log u}{x(u^x-1)}-\frac{\log(u^x-1)}{x^2}\Big)
$$ 
and $(\cdots)>0$ because $xu^x\log u=u^x\log(u^x)>(u^x-1)\log(u^x-1)$.
\kduk

\begin{defi}[$F_{m,n}(X)$]\label{def_Fmn}
Let $m,n\in\N$ with odd $n$. 
We define the function
$$
F_{m,n}(X)=\big((1+X)^m-X^m\big)^{1/n}\cc\R\to\R\,.
$$
\end{defi}

\begin{lemma}\label{lem_TaylIs}
Let $l,m,n\in\N$ and $l<m$. Then
$$
F_{m,\,n}(X)
=\sum_{j=0}^l\binom{m/n}{j}X^j+o(X^l)\ \ (X\to0)\,.
$$
\end{lemma}
\duk

\kduk

\begin{lemma}\label{lem_geneLege}
 Let $a\in\Z$, $m,d\in\N$, $p$ be a~prime with $(p,d)=1$ and let
 $$
 P_m(a,\,d):=\prod_{j=0}^{m-1}(a+jd)\,.
 $$
 Then
$$
\ord_p\big(P_m(a,\,d)\big)=
\sum_{j\ge1}\Big(\left\lfloor
\frac{m}{p^j}\right\rfloor+\ep_j\Big),\ \ep_j\in\{0,\,1\}\,.
$$
If $P_m(a,d)=m!$ then $\ep_j=0$ for every $j$.
\end{lemma}
\duk
Let $a$, $m$, $d$ and $p$ be as stated. For $j\in\N$ let $m_j$ 
($\in\N_0$) be the number of multiples of $p^j$ among the numbers $a+id$, 
$i=0,1,\ds,m-1$. A~double counting argument shows that
$$
\ord_p\big(P_m(a,\,d)\big)=\sum_{j\ge1}m_j\,.
$$
If integers $i$ and $i'$ are non-congruent modulo $p^j$, then so are 
$a+id$ and $a+i'd$. Hence if 
$I\sus\{0,1,\ds,m-1\}$ is an interval with length $|I|\le p^j$, then for at 
most one $i\in I$ the number $a+id$ is divisible by $p^j$, and if $|I|=p^j$ 
then there is  exactly one such number. Thus 
$$
m_j=\left\lfloor
\frac{m}{p^j}\right\rfloor+\ep_j
$$
with $\ep_j$ equal to $0$ or $1$, because we have the partition 
$$
\{0,\,1,\,\ds,\,m-1\}=I_1\cup I_2\cup\ds\cup I_k\cup I_0
$$
into intervals $I_1<I_2<\ds<I_k<I_0$ such that 
$k=\lfloor m/p^j\rfloor$, $|I_i|=p^j$
for $i>0$ and $|I_0|<p^j$. This proves the first claim. 
\kduk

\begin{lemma}\label{lem_denonBin}
Let $q$ be prime, $a\in\Z$ with $(q,a)=1$ and $k\in\N_0$. Then there exists $b\in\Z$ with $(q,b)=1$ such that
$$
\binom{a/q}{k}=\frac{b}{q^{k+\ord_q(k!)}}\,.
$$
\end{lemma}

\section[${}^c$The relation $q\,|\,x$]{The relation $q\,|\,x$}\label{sec_qdivx}

The following theorem was proven by Cassels in \cite{cass1}. We follow \cite[Chapter~6]{scho}.

\begin{thm}[Cassels, 1953]\label{thm_cass1}
If $p>q>2$ are primes and $x,y\in\Z^*$ are such numbers that 
$$
x^p-y^q=1\,,
$$ 
then $q$ divides $x$.    
\end{thm}
\duk
Suppose that $p$, $q$, $x$ and $y$ are as stated and that $(q,x)=1$. Then by Lemma~\ref{lem_oGCD} the
two factors in 
$$
(y+1)\cdot\frac{y^q+1}{y+1}=x^p
$$ 
are coprime. By PP1 (Proposition~\ref{prop_PP1} we have
$y+1=b^p$ for some $b\in\Z$. Since $x\ne0$, also $b\ne0$. Hence 
$$
x^p-(b^p-1)^q=1\,.
$$
We show that this equality cannot hold.

For $X\in\R$, we consider the function
$$
g(X)=X^p-(b^p-1)^q
$$
and show that $g(X)\ne1$ for every $X\in\Z$. Suppose that $b>0$. Since $y\ne0$, we have
$b\ge2$. Then
$$
g(b^q)=\sum_{j=0}^{q-1}b^{jp}(b^p-1)^{q-1-j}\ge q>1
$$
and
$$
g(b^q-1)=(b^q-1)^p-(b^p-1)^q<0\,,
$$
because, since $q<p$, by item~2 of Lemma~\ref{lem_incrDecr} we have 
$$
\big((b^q-1)^p\big)^{\frac{1}{pq}}=
(b^q-1)^{\frac{1}{q}}<
(b^p-1)^{\frac{1}{p}}=
\big((b^p-1)^q\big)^{\frac{1}{pq}}\,.
$$
The function $g(X)$ increases on $\R$ and we see that there is no $X\in\Z$ with $g(X)=1$.

Suppose that $b<0$, thus $b\le-1$. Then, similarly, 
$$
g(b^q)=\sum_{j=0}^{q-1}(b^p)^j(b^p-1)^{q-1-j}\ge q>1
$$
(each summand has sign $(-1)^{q-1}=1$) and
$$
g(b^q-1)=-((-b)^q+1)^p+((-b)^p+1)^q<0\,,
$$
because, since $q<p$, by item~1 of Lemma~\ref{lem_incrDecr} we have 
$$
\big(((-b)^q+1)^p\big)^{\frac{1}{pq}}=
((-b)^q+1)^{\frac{1}{q}}>
((-b)^p+1)^{\frac{1}{p}}=
\big(((-b)^p+1)^q\big)^{\frac{1}{pq}}\,.
$$
Again, $g(X)$ increases on $\R$ and we see that there is no $X\in\Z$ with 
$g(X)=1$. We reached a~contradiction and see that $q$ divides $x$.
\kduk

\section{The relation $p\,|\,y$}\label{sec_pdivy}

The following theorem was proven by Cassels in \cite{cass2}. We follow \cite[Chapter~6]{scho}.

\begin{thm}[Cassels, 1960]\label{thm_cass2}
If $p>q>2$ are primes and $x,y\in\Z^*$ are such numbers that 
$$
x^p-y^q=1\,,
$$ 
then $p$ divides $y$.    
\end{thm}
\duk
Let $p$, $q$, $x$ and $y$ be as stated and let $(p,y)=1$. By Lemma~\ref{lem_oGCD}, the
two factors in 
$$
(x-1)\cdot\frac{x^p-1}{x-1}=y^q
$$ 
are coprime. We deduce from this a~contradiction. 
By PP1 (Proposition~\ref{prop_PP1}),
$x-1=a^q$ with $a\in\Z$. Clearly $a\ne0$. Thus $y^q=(a^q+1)^p-1$ and with $F_{p,q}(X)$ as in Definition~\ref{def_Fmn} we express $y$ as
$$
y=a^p\cdot F_{p,q}(1/a^q)\,.
$$
We set $m=\lfloor p/q\rfloor+1$ ($\ge2$), $D=q^{m+\ord_q(m!)}$ and
$$
z=a^{mq-p}y-a^{mq}\cdot T_0^m(F_{p,q})(1/a^q)\ \ (\in\Q)\,.
$$ 
Using Lemmas~\ref{lem_TaylIs} and \ref{lem_denonBin} and the inequality $mq-p\ge0$ (following from $m>\frac{p}{q}$) we see that
$$
{\textstyle
Dz=Da^{mq-p}y-\sum_{k=0}^m
D\binom{p/q}{k}a^{mq-qk}\in\Z\,.
}
$$

We obtain a~contradiction by proving that the integer $Dz\ne0$ but at the same time $|Dz|<1$. Non-vanishing of $Dz$ follows from the non-divisibility 
$\neg(q\,|\,Dz)$: in the displayed expression for $Dz$ all terms are
divisible by $q$ except for the summand with $k=m$, which by
Lemma~\ref{lem_denonBin} is the
integer $D\binom{p/q}{m}$ not divisible by $q$.

We show that $|Dz|<1$. We have $z=a^{mq}\big(F_{p.q}(1/a^q)-
T_0^m(F_{p.q})(1/a^q)\big)$. Since $x\ne0$ but $q\,|\,x$ by
Theorem~\ref{thm_cass1}, $a\ne0,\pm1$ and $|a|\ge2$. We can use Lemma~\ref{lem_TaylBoun} with $X=1/a^q$ and get the bound
$$
|z|\le\frac{|a|^{mq}\cdot|a|^{-(m+1)q}}{(1-|a|^{-q})^2}=\frac{|a|^q}{(|a|^q-1)^2}\le\frac{1}{|a|^q-2}\le\frac{1}{|x|-3}\,.
$$
By Proposition~\ref{prop_lowBoux}, $|x |\ge q^{p-1}+q$. Thus 
$$
|Dz|\le q^{m+\ord_q(m!)-(p-1)}\,.
$$
If the exponent is negative, we are done. And indeed, by Lemma~\ref{lem_pOrdFac}, the inequality $m<\frac{p}{q}+1$ and since $p\ge5$ and $q\ge3$, it is negative:
\begin{eqnarray*}
m+\ord_q(m!)-(p-1)&\le&
{\textstyle
m\big(1+\frac{1}{q-1}\big)-(p-1)
}\\
&<&{\textstyle
\big(\frac{p}{q}+1\big)\big(1+\frac{1}{q-1}\big)-
(p-1)}\\
&=&{\textstyle
\frac{3-(p-2)(q-2)}{q-1}\le0\,.
}
\end{eqnarray*}
\kduk

\section{Corollaries}\label{sec_corrCass}

\chapter[Mih\u{a}ilescu's theorem: $x^p-y^q=1$ with $p>q>2$]{Mih\u ailescu's theorem, an outline of the proof}\label{chap_overview}

\section[${}^c$Theorems M0.95--M4]{Theorems M0.95--M4}

Catalan' conjecture claims that for integers $m,n\ge2$ the only 
nonzero solutions $x,y\in\Z^*$ of the Diophantine equation
$$
x^m-y^n=1
$$
occur for $m=2$ and $n=3$ as $\langle 
x,y\rangle=\langle\pm3,2\rangle$. In the first three chapters, we described 
the resolution of the conjecture in the case when $m=2$ or $n=2$. It is clear that to 
prove the whole conjecture, it remains to show that for every two distinct 
odd primes $p$ and $q$, the Diophantine equation
$$
x^p-y^q=1
$$
has no nonzero solution. This is what P.~Mih\u{a}ilescu accomplished in 
2004, and in the present chapter we outline his solution. We also describe the content of the
remaining chapters of our book. Recall that if $x,y\in\Z$ and $m,n\in\N$ are odd, then
$$
x^m-y^n=1\iff (-y)^n-(-x)^m=1\,,
$$
so that we may assume, if needed, that $p>q$.

In the first part of the chapter, we follow \cite[Chapter~1]{scho}, but we are more precise and clear about the assumptions in Theorems M1--M4. This pays off
because then we can simplify the deduction of Theorem~\ref{thm_mihai} from the four theorems compared to the argument given in \cite[p.~5]{scho}.

\begin{thm}[Mih\u{a}ilescu, 2004]\label{thm_mihai}
If $p$ and $q$ are distinct odd primes, then the Diophantine equation  
$$
x^p-y^q=1
$$
has no nonzero solution $x,y\in\Z^*$.
\end{thm}
Theorems~\ref{thm_mihai}, \ref{thm_x2y3} 
(or \ref{thm_eulerDesc}), \ref{thm_LebeThm} and \ref{thm_x2yq} 
together prove Catalan's conjecture.

Theorem~\ref{thm_mihai} follows from the next four Theorems M1--M4, all of 
which are due to P.~Mih\u{a}ilescu.

\begin{thm}[M1, \cite{miha1}]\label{thm_miha1}
If $p>q>2$ are primes such that $x^p-y^q=1$ for some numbers $x,y\in\Z^*$, then
$$
p^{q-1}\equiv1\ (\mathrm{mod}\;q^2)\,\text{ and }\,q^{p-1}\equiv1\ (\mathrm{mod}\;p^2)\,.
$$
\end{thm}
These congruences follow easily from the actual thing, which we dub as Theorem M0.95.

\begin{thm}[M0.95, \cite{miha1}]\label{thm_miha095}
If $p>q>2$ are primes such that $x^p-y^q=1$ for some numbers $x,y\in\Z^*$, then
$$
p^2\,|\,y\,\text{ and }\,q^2\,|\,x\,.
$$
\end{thm}

\begin{thm}[M2, \cite{miha}]\label{thm_miha2}
If $p,q\ge7$ are distinct primes such that $x^p-y^q=1$ for some numbers $x,y\in\Z^*$, then
$$
p\equiv1\ (\mathrm{mod}\;q)\,\text{ or }\,q\equiv1\ (\mathrm{mod}\;p)\,.
$$
\end{thm}
This is the hardest theorem of M1--M4  to prove. The original theorem in 
\cite{miha} holds for the primes $3$ and $5$ as well, but the present slightly 
restricted version has simpler proof and suffices to establish 
Theorem~\ref{thm_mihai}. 

\begin{thm}[M3, \cite{miha2}]\label{thm_miha3}
If $p,q\ge7$ are distinct primes such that $x^p-y^q=1$ for some numbers $x,y\in\Z^*$, then
$$
p<4q^2\,\text{ and }\,q<4p^2\,.
$$
\end{thm}
Again, the original theorem in 
\cite{miha2} holds for the primes $3$ and $5$ too.

\begin{thm}[M4, \cite{miha2}]\label{thm_M4}
If $p$ and $q$ are distinct odd primes and $p\in\{3,5\}$ or $q\in\{3,5\}$, then the Diophantine equation
$$
x^p-y^q=1
$$
has no nonzero solution $x,y\in\Z^*$.
\end{thm}
In \cite{miha2,scho} the theorem is proven for the larger set of primes $\{3,5,7,11,\ds,41\}$, but this restricted  version suffices for our purposes.

To deduce Theorem~\ref{thm_mihai} from Theorems M1--M4 we need a~simple lemma.

\begin{lemma}\label{lem_simpLema}
Let $q$ be prime, $x\in\Z$, $x\equiv1\ (\mathrm{mod}\;q)$ and $x^{q-1}\equiv1\ (\mathrm{mod}\;q^2)$. Then $x\equiv1\ (\mathrm{mod}\;q^2)$.    
\end{lemma}
\duk
We have $x=1+kq$ with $k\in\Z$. From 
$$
x^{q-1}=(1+kq)^{q-1}=1+\sum_{j=1}^{q-1}\binom{q-1}{j}(kq)^j\equiv1\ (\mathrm{mod}\;q^2)
$$
we deduce that $q^2\,|\,(q-1)kq$. Thus $q\,|\,k$, as stated. 
\kduk

\noindent
{\bf Deduction of Theorem~\ref{thm_mihai} from Theorems M1--M4. }In view of theorem M4 (Theorem~\ref{thm_M4}) we may assume that 
$p,q\ge7$ are distinct primes such that $x^p-y^q=1$ for some nonzero 
integers $x$ and $y$. We deduce a~contradiction. By Theorems M1 (Theorem~\ref{thm_miha1}), M2 (Theorem~\ref{thm_miha2}), Lemma~\ref{lem_simpLema} and the symmetry we have $p=1+kq^2$ for some $k\in\N$. By Theorem M3 (Theorem~\ref{thm_miha3}), 
$k\in\{1,2,3\}$. The values $k=1$ and $3$ are excluded because they give an 
even $p$. Thus $p=1+2q^2$. The values $q\ne3$ are excluded 
because they give a~$p$ divisible by $3$. We are left with the single pair $\langle p,q\rangle=\langle 19,3\rangle$. It is excluded by the assumption that $p,q\ge7$ (that is, by Theorem~\ref{thm_M4}). 
\kduk

\section{Overview of Chapters~\ref{chap_obstrGr}--\ref{chap_analysis}}

\chapter{An obstruction group}\label{chap_obstrGr}

\section{Number fields}\label{sec_NumFie}

A~field $L$ {\em extends} another field $K$, written $L/K$, if the base set of $K$ is a~subset of that of $L$, if both 
fields share neutral elements, which means that $0_K=0_L$ and $1_K=1_L$, and if  
addition and multiplication in $L$ extend these operations in~$K$. 
In this situation $L$ is a~$K$-vector space. We denote its dimension by 
$[L:K]$.

\begin{defi}[number fields]\label{def_numFie}
A~field extension $K$ of the field of fraction $\Q_{\mathrm{fi}}$ with 
finite dimension $d:=[K:\Q_{\mathrm{fi}}]$ ($\in\N$) is called a~number field (with degree $d$).
\end{defi}

Let $K$ be a~number field with degree $d$ and let $\al\in K$. Then the $d+1$ 
powers $1_K$, $\al$, $\al^2$, $\ds$, $\al^d$ are linearly dependent over 
$\Q$, and $\al$ is a~root of a~nonzero polynomial in $\Q[x]$ of degree at 
most $d$.  

\begin{defi}[rings of integers]\label{def_rinInt}
Let $K$ be a~number field. We define the subset
$$
\mathcal{O}_K=\{\al\in K\cc\;\text{$p(\al)=0_K$ for a~monic polynomial $p(x)\in\Z[x]$}\}\,.
$$
We say that $\mathcal{O}_K$ is the (base set of the) ring of integers of $K$.
\end{defi}

\begin{prop}\label{prop_OK}
The set $\mathcal{O}_K$ contains the neutral elements $0_K$ and $1_K$, and is closed under addition and multiplication in $K$. It forms a~subdomain of $K$.    
\end{prop}
\duk

\kduk

\begin{prop}\label{prop_h1AndUFD}
For every number field $K$ its ring of integers is {\em UFD} if and only if the class number $h_K=1$.     
\end{prop}
\duk

\kduk

\section{Cyclotomic fields}\label{sec_CyclFie}

\begin{prop}\label{prop_cyclInteg}
For every prime $p$ the ring of integers of the cyclotomic field $\Q(\zeta_p)_{\mathrm{fi}}$ equals $\Z[\zeta_p]_{\mathrm{do}}$.    
\end{prop}
\duk

\kduk

\begin{thm}\label{thm_z3Az5}
The domains
$$
\Z[\zeta_3]_{\mathrm{do}}\,\text{ and }\,\Z[\zeta_5]_{\mathrm{do}}
$$
are Euclidean. In view of Propositions~\ref{prop_h1AndUFD} and \ref{prop_cyclInteg}, and of 
Theorem~\ref{thm_euclIsUFD} the class numbers of cyclotomic fields $\Q(\zeta_3)_{\mathrm{fi}}$ and $\Q(\zeta_5)_{\mathrm{fi}}$ are
$$
h_3=h_5=1\,.
$$
\end{thm}
\duk

\kduk

\section{The obstruction group}\label{sec_obsGroup}

\chapter[Super--Cassels relations: $p^2\,|\,y$ and $q^2\,|\,x$]{Super--Cassels relations}\label{chapsupercass}

\section{The Stickelberger ideal}

\section{Theorems M$0.95$ and  M$1$}

\chapter[Theorem M4: $p=3,5$ or $q=3,5$]{Theorem M4}\label{chap_M4}

\begin{thm}[M4]\label{thm_miha4}
If $p$ and $q$ are distinct odd primes and $p\in\{3,5\}$ or $q\in\{3,5\}$, then the Diophantine equation
$$
x^p-y^q=1
$$
has no nonzero solution $x,y\in\Z^*$.
\end{thm}
\duk

\kduk

\appendix
\chapter{Results from mathematical analysis}\label{chap_analysis}

\section{Taylor series}

\begin{thm}[on Taylor polynomials]\label{thm_TaylorPol}
Let $n\in\N$, $I\sus\R$ be an open interval, $a\in I$ and $f\cc I\to\R$ be a~function that has $n$ derivatives
$$
f',\,f'',\,\ds,\,f^{(n)}\cc I\to\R\,.
$$
Then the (Taylor) polynomial
$$
T_f^{n,\,a}(x)
:=\sum_{j=0}^n\frac{f^{(j)}(a)}{j!}(x-a)^j
$$
is a~unique polynomial $P\in\R[x]$ with degree at most $n$ and the 
property that 
$$
f(x)=P(x)+o((x-a)^n)\ \ (x\to a)\,.
$$
\end{thm}
\duk

\kduk

\chapter{Results from algebra}\label{chap_algebra}

\section{Commutative group rings}

By a~group we mean an Abelian and multiplicatively written group
$$
G=\langle G,\,1_G,\,\cdot\rangle
$$
and by a~ring we mean a~commutative unital ring 
$$
R=\langle R,\,0_R,\,1_R,\,+,\,\cdot\rangle\,. 
$$

\begin{defi}[group rings]
Let $R$ be a~ring and $G$ be a~group. The group ring $R[G]$ generated by $R$ and $G$ is the algebraic structure
$$
R[G]=\langle S,\,0_S,\,1_S,\,+,\,\cdot\rangle\,,
$$
with the base set, neutral elements, and operations defined as follows. $S$ 
is the set of maps $F\cc G\to R$ such that $F(g)\ne 0_R$ for only finitely 
many $g\in G$. $0_S$ is the zero map with $0_S(g)=0_R$ for every $g\in G$. 
$1_S$ is the map with $1_S(1_G)=1_R$ and $1_S(g)=0_R$ for every $g\in 
G\setminus\{1_G\}$. The addition is defined component-wise by
$$
(F+F')(g)=F(g)+F'(g),\ g\in G\,.
$$
Finally, the multiplication is defined 
by the convolution
$$
(F\cdot F')(g)=\sum_{\substack{e,\,f\in G\\e\cdot f=g}}F(e)\cdot F'(f),\ g\in G\,,
$$
where the possibly infinite sum is defined as the sum, in $R$, of 
the finitely many summands different from $0_R$. 
\end{defi}

\begin{thm}
$R[G]$ is a~commutative unital 
ring.     
\end{thm}
\duk

\kduk


\addcontentsline{toc}{chapter}{References}
\begin{thebibliography}{100}

\bibitem{bilu_al}
Yu. Bilu, Y. Bugeaud and M. Mignotte, {\em The Problem of Catalan},
Springer, Cham, 2014

\bibitem{casselsWiki}
J. W. S. Cassels, Wikipedia article, \url{https://en.wikipedia.org/wiki/J._W._S._Cassels}

\bibitem{cass1}
J.\,W.\,S. Cassels, On the equation $a^x-b^y=1$, {\em Amer. J. Math.} {\bf 75} (1953), 159--162

\bibitem{cass2}
J.\,W.\,S. Cassels, On the equation $a^x-b^y=1$. II, {\em Proc. Cambridge Phil. Soc.} {\bf 56} (1960), 97--103; Corrigendum: Ibid, {\bf 57} (1961), 187 

\bibitem{cata1}
E. Catalan, Probl\`eme 48, {\em Nouv. Ann. Math.} {\bf 1} (1842), 520

\bibitem{cata2}
E. Catalan, Note extraite d'une lettre adress\'ee \`a l'\'editeur, {\em J. reine angew. Math.} {\bf 27} (1844), 192

\bibitem{chein}
E.\,Z. Chein, A~note on the equation $x^2=y^q+1$, {\em Proc. Amer. Math. Soc.} {\bf 56} (1976), 83--84

\bibitem{cohe1}
H. Cohen, {\em Number Theory. Volume I: Tools and Diophantine Equations}, Springer, New York 2007 

\bibitem{conr}
B. Conrad, An example of descent by Euler, 5 pp., \url{https://kconrad.math.uconn.edu/blurbs/ugradnumthy/descentbyeuler.pdf}\\ (visited in December 4, 2025)

\bibitem{EKR}
P. Erd\H{o}s, Chao Ko and R. Rado, 
Intersection theorems for systems of finite sets,
{\em Quart. J. Math. Oxford}. Ser. (2) {\bf 12} (1961), 313--320

\bibitem{EKRwiki}
Erd\H{o}s--Ko--Rado theorem, Wikipedia article, \url{https://en.wikipedia.org/wiki/Erd%C5%91s%E2%80%93Ko%E2%80%93Rado_theorem}

\bibitem{EucDom_wiki}
Euclidean domain, Wikipedia article, \url{https://en.wikipedia.org/wiki/Euclidean_domain}

\bibitem{cataWiki}
Eug\`ene Charles Catalan, Wikipedia article,
\url{https://en.wikipedia.org/wiki/Eug%C3%A8ne_Charles_Catalan}

\bibitem{eulerWiki}
Leonhard Euler, Wikipedia article, \url{https://en.wikipedia.org/wiki/Leonhard_Euler}

\bibitem{eule}
L. Euler, Theorematum quorundam arithmeticorum demonstrationes, {\em Comm. Acad. Sci. Petrop.} {\bf 10} (1738), 125--146

\bibitem{hard_wrig}
G.\,H. Hardy and E.\,M. Wright, {\em An Introduction to the Theory of
Numbers}. Fourth Edition, Oxford at the Clarendon Press, Oxford 1960 

\bibitem{chaoKoWiki}
Ke Zhao, Wikipedia article, \url{https://en.wikipedia.org/wiki/Ke_Zhao}

\bibitem{klaz}
M. Klazar, O \v re\v sen\'\i\ diofantick\'e 
rovnice $x^2-y^3=\pm1$, {\em Matematick\'e  obzory} {\bf 32} (1989), 47--53 (Solving the Diophantine equation $x^2-y^3=\pm1$)

\bibitem{ko}
Chao Ko, On the Diophantine equation $x^2=y^n+1$, $xy\ne0$, {\em Sci. Sinica} {\bf 14} (1965), 457--460

\bibitem{lang_algebra}
S. Lang, {\em Algebra}. Revised Third Edition, Springer-Verlag, New York 2002 

\bibitem{lebe}
V. Lebesgue, Sur l'impossibilit\'e en nombres entiers de l'\'equation $x^m=y^2+1$, {\em  Nouv. Ann. Math.} {\bf 8} (1850), 178--181

\bibitem{Vlebe}
Victor Am\'ed\'ee Lebesgue, Mac Tutor article,
\url{https://mathshistory.st-andrews.ac.uk/Biographies/Lebesgue_Victor/}

\bibitem{mets}
T. Mets\"ankyl\"a, Catalan's conjecture: another old Diophantine problem solved, {\em Bulletin Amer. Math. Soc.} {\bf 41} (2003), 43--57

\bibitem{miha1}
P. Mih\u{a}ilescu, A~class number free criterion for Catalan's conjecture, {\em 
J.~Number Theory} {\bf 99} (2003), 225--231 

\bibitem{miha}
P. Mih\u{a}ilescu, Primary cyclotomic units and a~proof of Catalan's 
conjecture, {\em J.~reine angew. Math.}  {\bf 572} (2004), 167--195

\bibitem{miha2}
P. Mih\u{a}ilescu, On the class group of cyclotomic extensions in presence 
of a~solution to Catalan's equation, {\em J. Number Theory} {\bf 118} 
(2006), 123--144

\bibitem{mord}
L.\,J. Mordell, {\em Diophantine Equations}, Academic Press, London and New York 1969

\bibitem{nage}
T. Nagell, Sur l'impossibilit\'e de l'\'equation ind\'etermin\'ee
$z^p+1=y^2$, {\em Norsk. Mat. Forenings Skrifter} {\bf 4} (1921),
14 pp.

\bibitem{nota}
Ch. Notari, Une r\'esolution \'el\'ementaire de l'\'equation diophantienne $x^3=y^2-1$, {\em Expositiones Mathematicae} {\bf 21} (2003), 279--283

\bibitem{ribe}
P. Ribenboim, {\em Catalan's Conjecture.
Are 8 and 9 the Only Consecutive Powers?}, 
Academic Press, Inc., Boston, MA, 1994

\bibitem{roge}
K. Rogers, The axioms for Euclidean domains, {\em Amer. Math. Monthly} {\bf 78} (1971), 1127--1128

\bibitem{scho}
R. Schoof, {\em Catalan's Conjecture}, 
Springer-Verlag London, Ltd., London, 2008

\bibitem{wakuliczWiki}
Antoni Wakulicz, Wikipedia article, \url{https://pl.wikipedia.org/wiki/Antoni_Wakulicz} (in Polish)

\bibitem{waku}
A. Wakulicz, On the equation $x^3+y^3=2z^3$, {\em Colloq. Math.} {\bf 5} (1957), 11--15 

\bibitem{weil}
A. Weil, {\em Number Theory. An Approach through History: From Hammurapi to Legendre}, Birhkh\"auser, Boston 1984

\end{thebibliography}
\end{document}